\documentclass[final]{siamltex}

\usepackage{amssymb}
\usepackage{graphicx}
\usepackage{mathrsfs}

\def\3n{\negthinspace \negthinspace \negthinspace }
\def\2n{\negthinspace \negthinspace }
\def\1n{\negthinspace }

\def\ds{\displaystyle}

\def\={\buildrel \triangle \over =}

\def\a{\alpha}

\def\e{\varepsilon}

\def\o{\omega}

\def\i{\infty}

\def\L{\Lambda}

\def\O{\Omega}

\def\cA{\mathscr{A}}

\def\cF{\mathscr{F}}

\def\cM{\mathscr{M}}
\def\cN{\mathscr{N}}

\def\cS{\mathscr{S}}
\def\cT{\mathscr{T}}

\def\ms{\medskip}

\def\q{\quad}
\def\qq{\qquad}

\def\ra{\mathop{\rightarrow}}

\def\esssup{\mathop{\rm esssup}}
\def\essinf{\mathop{\rm essinf}}
\def\max{\mathop{\rm max}}
\def\min{\mathop{\rm min}}

\def\sup{\mathop{\rm sup}}

\def\tilde{\widetilde}
\def\hat{\widehat}
\def\wt{\widetilde}

\def\|{\Big |}
\def\({\Big (}
\def\){\Big )}
\def\[{\Big[}
\def\]{\Big]}

\newtheorem{remark}[theorem]{Remark}

\def\be{\begin{equation}}
\def\bea{\begin{eqnarray*}}
\def\eea{\end{eqnarray*}}
\def\bel{\begin{equation}\label}
\def\ee{\end{equation}}
\def\bt{\begin{theorem}}
\def\bcd{\begin{condition}}
\def\ecd{\end{condition}}
\def\et{\end{theorem}}
\def\bc{\begin{corollary}}
\def\ec{\end{corollary}}
\def\bde{\begin{definition}}
\def\ede{\end{definition}}
\def\bl{\begin{lemma}}
\def\el{\end{lemma}}
\def\bp{\begin{proposition}}
\def\ep{\end{proposition}}
\def\br{\begin{remark}}
\def\er{\end{remark}}
\def\ba{\begin{array}}
\def\ea{\end{array}}
\def\ed{\end{document}}

\def\ds{\displaystyle}

\title{Optimal Switching of One-Dimensional Reflected BSDEs, and
Associated Multi-Dimensional BSDEs with Oblique Reflection}

\author{Shanjian Tang\thanks{Department of
Finance and Control Science, School of Mathematical Sciences, Fudan
University, Shanghai 200433, China; and Graduate Department of Financial
Engineering, Ajou University, San 5, Woncheon-dong, Yeongtong-gu, Suwon,
443-749, Korea. This author is supported in part by NSFC Grant 10325101,
Basic Research Program of China (973 Program) Grant 2007CB814904, by the
Science Foundation of the Ministry of Education of China Grant
\#200900071110001,by Chang Jiang Scholars Programme, and by WCU (World Class
University) Program through the Korea Science and Engineering Foundation
funded by the Ministry of Education, Science and Technology (R31-20007) ({\tt
sjtang@fudan.edu.cn})} \and Wei Zhong\thanks{Institute of Mathematics and
Department of Finance and Control Science, School of Mathematical Sciences,
Fudan University, Shanghai 200433, China ({\tt zhongwei@fudan.edu.cn})} \and
Hyeng Keun Koo\thanks{Graduate Department of Financial Engineering, Ajou
University, San 5, Woncheon-dong, Yeongtong-gu, Suwon, 443-749, Korea.  This
author is supported by WCU (World Class University) Program through the Korea
Science and Engineering Foundation funded by the Ministry of Education,
Science and Technology (R31-20007) ({\tt hkoo@ajou.ac.kr})}}

\begin{document}

\maketitle

\begin{abstract}
In this paper, an optimal switching problem  is proposed for one-dimensional
reflected backward stochastic differential equations (RBSDEs, for short)
where  the generators, the terminal values and the barriers are all switched
with positive costs. The value process is characterized by a system of
multi-dimensional RBSDEs with oblique reflection, whose existence and
uniqueness are by no means trivial and are therefore carefully examined.
Existence is shown using both methods of the Picard iteration and
penalization, but under some different conditions. Uniqueness is proved by
representation either as the value process to our optimal switching problem
for one-dimensional RBSDEs, or as the equilibrium value process to a
stochastic differential game of switching and stopping. Finally, the switched
RBSDE is interpreted as a real option.
\end{abstract}



\pagestyle{myheadings} \thispagestyle{plain} \markboth{S. TANG, W.
ZHONG AND H. K. KOO}{OPTIMAL SWITCHING OF REFLECTED BSDES}

\section{Introduction}

Let $\{W(t), 0\leq t\leq T\}$ be a $d$-dimensional standard Brownian
motion defined on some complete probability space $(\O,\cF,P)$, and
denote by $\{\cF_t,\, 0\leq t\leq T\}$ the natural filtration,
augmented by all the $P$-null sets of $\cF$. Define the following
spaces of real-valued processes. \begin{eqnarray*}\ba{rcl}
\cS'^2&\triangleq&\displaystyle  \left\{ \phi: \phi \hbox{ \rm is
}\{\cF_t, 0\le t\le T\}\mbox{\rm -adapted and r.c.l.l. }
s.t.\, E[\sup\limits_{0\leq t\leq T}|\phi(t)|^2]<\infty\right\},\\
\cS^2&\triangleq &\{ \phi\in\cS'^2:
\phi\mbox{ is continuous}\},\\
\cN'^{\, 2}&\triangleq &\{\phi\in \cS'^2: \phi\mbox{ is increasing and
}\phi(0)=0\},\\
\cN^{\,2}&\triangleq &\{\phi\in \cN'^{\,2}: \phi\mbox{ is
continuous}\},\\\cM^2&\triangleq&\{\phi: \phi \hbox{ \rm is }
\{\cF_t, 0\le t\le T\}\mbox{\rm -predictable and square-integrable
on $[0,T]\times \Omega$} \}.\ea\end{eqnarray*}

Let $\{\theta_j\}_{j=0}^{\i}$ be an increasing sequence of stopping times
with values in $[0,T]$. For any $j,\, \alpha_j$ is an
$\cF_{\theta_j}$-measurable random variable with values in
$\L\triangleq\{1,\cdots,m\}$. Assume that $a.s. \,\o$, there exists an
integer $N(\o)<\i$ such that $\theta_N=T$. Then we define an admissible
switching  as:
$$a(s)=\alpha_0\chi_{[\theta_0,\,\theta_1]}(s)+\sum_{j=1}^{N-1}\alpha_j\chi_{(\theta_j,\, \theta_{j+1}]}(s),\q s\in[\theta_0,\,T].$$
Denote  by $\cA_t^i$ all the admissible switching control with initial data
$\alpha_0=i\in\L,\,\theta_0=t$. For given $a\in\cA_t^i,\,
\xi=(\xi_1,\cdots,\xi_m)^T\in L^2(\O, \cF_T, P; R^m)$ and
$S=(S_1,\cdots,S_m)^T\in (\cS^2)^m$, consider the following switched
reflected backward stochastic differential equation (abbreviated as RBSDE):
\be\label{BSDE1}\left\{\ba{rcl} U^a(s)&=&\displaystyle
\xi_{a(T)}+\int_s^T\psi(r,U^a(r),V^a(r),a(r))d
r- (L^a(T)-L^a(s))\\[3mm]&&\displaystyle-\sum\limits_{j=1}^{N-1}
 [U^a(\theta_j)-h_{\alpha_{j-1},\,\alpha_j}(\theta_j,U^a(\theta_j))]\chi_{(s,T]}(\theta_j)\\[3mm]
 &&\displaystyle
-\int_s^T V^a(r)d
W(r), \qq s\in[t,T];\\[3mm]U^a(s)&\leq & S_{a(s)}(s), \qq s\in[t,T]; \\[3mm]&&\displaystyle\int_t^T(U^a(s)-S_{a(s)}(s))\, d L^a(s)=0. \ea\right.\ee
Here and in the following, $\chi$ is an indicator function. The generator
$\psi$, the terminal condition $\xi$ and the upper barrier $S$ of RBSDE
(\ref{BSDE1}) are all switched by $a$. At each switching time $\theta_j$
before termination, the value of $U^a$ will jump by an amount of
$U^a(\theta_j)-h_{\alpha_{j-1},\,\alpha_j}(\theta_j,U^a(\theta_j))$ which can
be regarded as a penalty or cost for the switching. RBSDE (\ref{BSDE1})  can
be solved in a backwardly recursive way in the subintervals
$[\theta_{N-1},T]$ and $[\theta_{j-1}, \theta_j)$ for $j= N-1,\cdots,2,1$. To
be precise, RBSDE~ (\ref{BSDE1}) in the last subinterval $[\theta_{N-1},T]$
reads: \be\label{BSDE1.1}\left\{\ba{rcl} U^a(s)&=&\displaystyle
\xi_{a(T)}+\int_s^T\psi(r,U^a(r),V^a(r),\alpha_{N-1})\, d r
\\[3mm]&&\displaystyle-(L^a(T)-
L^a(s))-\int_s^T V^a(r)\, d
W(r), \qq s\in[\theta_{N-1},T];\\[3mm]U^a(s)&\leq & S_{\a_{N-1}}(s), \qq s\in[\theta_{N-1},T];
\\[3mm]&&\displaystyle\int_{\theta_{N-1}}^T(U^a(s)-S_{\a_{N-1}}(s))\, d L^a(s)=0.\ea\right.\ee
From~\cite[Theorem 5.2]{El Karoui}, RBSDE~(\ref{BSDE1.1}) has a unique
adapted solution on $[\theta_{N-1}, T]$ under Hypothesis 1 (see Section 2
below). Then  we have
$$U^a(\theta_{N-1}-)=h_{\a_{N-2}, \a_{N-1}}(\theta_{N-1}, U^a(\theta_{N-1})),$$
which serves as the terminal value of RBSDE~ (\ref{BSDE1}) in $[\theta_{N-2},
\theta_{N-1})$. In general, in the subinterval $[\theta_{j-1}, \theta_j)$ for
$j= N-1,\cdots,2,1$, RBSDE~(\ref{BSDE1}) reads
\be\label{BSDE1.2}\left\{\ba{rcl} U^a(s)&=&\displaystyle
h_{\alpha_{j-1},\,\alpha_j}(\theta_j,U^a(\theta_j))+\int_s^{\theta_j}\psi(r,U^a(r),V^a(r),\alpha_{j-1})\,
d
r\\[3mm]&&\displaystyle-(L^a(\theta_j)-L^a(s))-\int_s^{\theta_j} V^a(r)\, d
W(r),\qq
s\in[\theta_{j-1},\theta_j);\\[3mm]U^a(s)&\leq & S_{\a_{j-1}}(s),\qq
s\in[\theta_{j-1},\theta_j);\\[3mm]&&\displaystyle\int_{\theta_{j-1}}^{\theta_j}(U^a(s)-S_{\a_{j-1}}(s))\, d L^a(s)=0.\ea\right.\ee
 Here $U^a(\theta_j)$ is
specified in the interval $[\theta_j, \theta_{j+1})$ and we have the
following relations under Hypothesis 2 (see Section 2 below):
$$h_{\alpha_{j-1},\,\alpha_j}(\theta_j,U^a(\theta_j))\leq
U^a(\theta_j)\leq S_{\a_{j-1}}(\theta_j).$$  The existence and uniqueness of
an adapted solution to RBSDE~(\ref{BSDE1}) in the interval $[0, T]$ are
obtained in an obvious way from the existence and uniqueness of an adapted
solution to RBSDE~(\ref{BSDE1}) in all the subintervals $[\theta_{N-1},T]$
and $[\theta_{j-1}, \theta_j)$ for $j= N-1,\cdots,2,1$.

One-dimensional RBSDEs were introduced  with  fixed single reflecting
barrier---as a generalization to the associated right Hamilton-Jacobi-Bellman
equation for the value process of a traditional optimal stopping problem for
SDEs---first by El Karoui et al.~\cite{El Karoui}, who gave the first
existence and uniqueness results for one-dimensional RBSDEs. They were later
generalized by Cvitanic and Karatzas~\cite{Cvitanic} to the case of fixed
double reflecting barriers and linked  to the well-known Dynkin games. Recent
literature shows an interest in using a one-dimensional RBSDE to specify the
cost functional in the modelling of optimal stochastic control/stochastic
differential games. See Wu and Yu~\cite{WuYu}, who discussed the {\it
classical} (by which we mean that any instant action of control could cause a
jump neither to the system state nor to the cost) optimal stochastic control
problem for a system of SDEs and a one-dimensional RBSDE.

In this paper, we study the general optimal switching ({\it non-classical})
problem for one-dimensional RBSDE~(\ref{BSDE1}), where the generator, the
terminal value and the upper barrier are all switched with positive costs.
Such a model appears in the management of real options. A real option is the
right, but not the obligation, to undertake some business decisions;
typically the option to make, abandon, expand, or contract a capital
investment. The first solution component $U^a$ of RBSDE~(\ref{BSDE1}) can be
interpreted as the minimal value process of a real option---also as the
wealth process of one counterparty (called the holder hereafter)---which is
subject to a discretionary calling back or termination by the other
counterparty (called the issuer) like at any time $\tau$ at a cost of paying
the amount of money $S(\tau)\chi_{\tau<T}+\xi\chi_{\tau=T}$ to the holder:
the set $\L$ models the totality of instant possible decision choices, and
the switching $a$ represents the holder's underlying decision process. Our
optimal switching problem is then the optimal management of the real option
for the holder to maximize the above minimal value process of the real option
by suitably and dynamically making his/her decisions. We will discuss real
options related to our switched RBSDEs in more detail in Section 6.

Mathematically,  the optimal switching problem for RBSDE~(\ref{BSDE1}) is to
maximize $U^a(t)$ over $a\in \cA_t^i, i\in \L$. The value process turns out
to be described by the following system of multi-dimensional RBSDEs
 with double reflecting barriers: for $i\in
\L$, \be\label{DRBSDE}\left\{\ba{rcl}
Y_i(t)&=&\displaystyle\xi_i+\int_t^T\psi(s,Y_i(s),Z_i(s),i)d s-\int_t^T d
K_i^-(s)\\[3mm]&&\displaystyle+\int_t^T d K_i^+(s)-\int_t^T Z_i(s)d W(s),\q t\in [0,T];\\[3mm]
S_i(t)&\geq&Y_i(t)\geq \max\limits_{j\neq i, j\in \L}h_{i,j}(t,Y_j(t)),\q t\in [0,T];\\[3mm]
&&\displaystyle\int_0^T\left(Y_i(s)-\max\limits_{j\neq i, j\in
\L}h_{i,j}(s,Y_j(s))\right)\, dK_i^+(s)=0,\\[3mm]
&&\displaystyle\int_0^T\left(Y_i(s)-S_i(s)\displaystyle\right)\, d
K_i^-(s)=0.\ea\right.\ee The last two equalities are respectively called the
{\it lower} and the {\it upper minimal boundary conditions}. Solution of the
above RBSDE~(\ref{DRBSDE}) is by no means trivial, and will be examined
carefully in this paper. The unusual feature here is that for
RBSDE~(\ref{DRBSDE}), the upper barrier is a fixed process, while the lower
barrier depends on the unknown process and is therefore implicit, which is
different from one-dimensional RBSDEs with fixed double barriers. In contrast
to RBSDEs with oblique reflection introduced in Hu and Tang~\cite{Hu2}, there
is an additional fixed upper barrier. This difference will complicate the
analysis of the existence and uniqueness for solutions to
RBSDE~(\ref{DRBSDE}). For $t\in[0,T]$, define
$$Q(t)\triangleq\{(y_1,\cdots, y_m)^T\in R^m : h_{i,j}(t,y_j) \leq
y_i\leq S_i(t) ,\forall \,i,j\in \L, j\neq i\}.$$ Then the state process
$Y(\cdot)$ of (\ref{DRBSDE}) is forced to evolve in the time-dependent set
$Q(\cdot)$, thanks to the accumulative action of two increasing processes
$K_i^+$ and $K_i^-$.

The literature on RBSDEs exhibits an interest in both methods of the
penalization and the the Picard iteration--- like in the pioneering work on
one-dimensional RBSDE of ~El Karoui et al.~\cite{El Karoui}, and in the very
recent works on optimal switching of one-dimensional BSDEs by Hu and
Tang~\cite{Hu2} for the penalty method and by Hamad\'ene and
Zhang~\cite{Hamadene4} for the Picard iteration method. This motivates us to
develop (in Sections 3 and 4, respectively) both methods to the existence  of
the solution of RBSDE~(\ref{DRBSDE})--- which is the associated right
Hamilton-Jacobi-Bellman equation for the the optimal switching of
one-dimensional RBSDEs. Assuming that the fixed barrier is super-regular, we
firstly prove existence of the solution by the method of a Picard
 iteration, invoking a generalized monotonic limit theorem.
 As a key condition of the generalized monotonic limit theorem,
 the comparison of the increment of the increasing process is necessary,
 which is formulated as Lemma~\ref{l1}.
 In addition, the proof of the minimal boundary condition is complicated by the appearance of the additional fixed barrier,
 and our arguments look very technical and seem to be new. Secondly we consider the case of a particular  barrier,
 and obtain an existence result by a penalty method, without assuming the super-regularity on the fixed barrier,
 which on one hand exhibits an advantage over the former method of Picard iteration.  On the
 other hand, our penalty method is not able to
 treat as general switching costs as the method of Picard iteration deals with, which will help to stimulate further efforts
 to decrease or even to remove such a disadvantage.

 It is also worth noting that though it itself
 looks shorter, the presentation of the Picard iteration method invokes the generalized monotonic limit theorem~\cite[Theorem 3.1]{Peng} (see Lemma~\ref{T2.1} in Section 2 below), a
 comprehensive proof of which is very long---in fact, it has to further employ very technical lemmas of Peng~\cite{Peng0}.

 Uniqueness of the solution to RBSDE~(\ref{DRBSDE}) is proved in Section 5 by linking
 it either to the value process for our optimal switching of one-dimensional RBSDEs,
 or to a stochastic game involving both switching and stopping control for one-dimensional BSDEs.


 The rest of the paper is
 organized as follows: In Section 2,
 we formulate our problem, introduce the generalized monotonic
 theorem,
 and give some preliminary results on RBSDEs, which will be used in subsequent arguments.
 In Section 3, we prove existence of the solution by the method of Picard iteration.
 In Section 4, existence of the solution is shown by the
 penalty method. Uniqueness of the
 solution is shown in Section 5. We discuss economic
interpretations in Section 6 with an example.

\section{Preliminaries}

We make the following assumption on the generator
$\{\psi(\cdot,\cdot,\cdot,i),\\ i\in\L\}$.

\ms  \textbf{Hypothesis 1. } The generator $\psi$ satisfies the following:

(i) The process $\psi(\cdot,y,z,i)\in \cM^2 $ for any $ (y,z, i)\in R\times
R^d\times \L$.

(ii) There is a constant $C>0$ such that for $(y,y',z,z')\in R\times R\times
R^d\times R^d$ and  $(t,i)\in [0,T]\times \L$, we have  $$\q
|\psi(t,y,z,i)-\psi(t,y',z',i)|\leq C(|y-y'|+|z-z'|).$$

\ms We make the following two assumptions on the function $\{h_{i,j},
i,j\in\L\}$, which are introduced in ~\cite{Hamadene4}.

\ms \textbf{Hypothesis 2. } For any $(i,j)\in \L\times \L,$ the function
$h_{i,j}(t,y)$ is continuous in $(t,y)$, increasing in $y$, and
$h_{i,j}(t,y)\leq y.$

\ms \textbf{Hypothesis 3. } For any $y_n\in R$ and any loop $\{i_k\in\L,\,
k=1,\cdots,n\}$ such that $i_1=i_n$ and $i_k\neq i_{k+1}$ for
$k=1,2,\ldots,n-1$, define $y_k\triangleq h_{i_k,i_{k+1}}(t,y_{k+1})$ for
$k=1,\cdots,n-1.$ Then $y_1<y_n$.

\ms In Section 4, we shall specialize the  function $h_{i,j}$ to the form:
$h_{i,j}(t,y)=y-k(i,j)$ for some positively valued function $k$ defined on
$\L\times \L$. We shall make the following assumption on $k$, which is
introduced in ~\cite{Hu2}.

\ms \textbf{Hypothesis 3'. } The function $k: \L\times \L\to R$ satisfies the
following:

(i) $\forall \,(i,j)\in\L\times \L, \q k(i,j)>0$ for $i\neq j,$ and $
k(i,i)=0$.

(ii) $\forall \,(i,j,l)\in\L\times\L\times\L$ such that $i\neq j,\, j\neq l,
\q k(i,j)+k(j,l)\geq k(i,l).$

\ms  \br Hypothesis 3 means that there is no free loop of instantaneous
switchings. Hypotheses 2 and 3 are satisfied when $h_{i,j}(t,y)=y-k(i,j)$ for
$(t,y)\in [0,T]\times R$ and $(i,j)\in \L\times\L$ with the function $k$
satisfying Hypothesis 3' (i).\er

\bde An adapted solution of system (\ref{DRBSDE}) is a quadruple
$$\ba{rcl}(Y,Z,K^+,K^-) &\triangleq& \{Y(t),Z(t), K^+(t),K^-(t);\, 0\leq t\leq
T\}\\
&\in& (\cS^2)^m\times (\cM^2)^{m\times d}\times (\cN^{\,2})^m\times
(\cN^{\,2})^m,\ea $$  taking values in $R^m\times R^{m\times d}\times
R^m\times R^m$ and satisfying~(\ref{DRBSDE}). \ede

We recall here the generalized monotonic limit theorem \cite[Theorem
3.1]{Peng}, which will be used in our method of Picard iteration.

\bl\label{T2.1}\textbf{(Generalized monotonic theorem)} We assume
the following sequence of It\^{o} processes: \begin{eqnarray*}
y^n(t)=y^n(0)+\int_0^t g^n(s)\, ds-K^{+,n}(t)+K^{-,n}(t)+\int_0^t
z^n(s)\, dW(s),\q n=1,2,\cdots\end{eqnarray*} satisfy

(i) for each $n, \, g^n\in \cM^2, K^{+,n}\in \cN^{\,2}, K^{-,n}\in
\cN'^{\,2}$;

(ii) $K^{-,n_2}(t)-K^{-,n_2}(s)\geq K^{-,n_1}(t)-K^{-,n_1}(s), \forall\,
0\leq s\leq t\leq T, \forall \,n_1\leq n_2$;

(iii) For each $t\in[0,T],\,  \{K^{-,n}(t)\}_{n=1}^\i$ increasingly converges
to $K^-(t)$ with
$$E|K^-(T)|^2<\i;$$

(iv) $(g^n,z^n)_{n=1}^\i$ converges to $(g, z)$ weakly in $\cM^2$;

(v)  For each $t\in[0,T],\, \{y^n(t)\}_{n=1}^\i$ increasingly converges  to
$y(t)$ with
$$E\sup\limits_{0\leq t\leq
T}|y(t)|^2<\infty.$$

Then the limit $y$ of $\{y^n\}_{n=1}^\i$ has the following form
\begin{eqnarray*} y(t)=y(0)+\int_0^t g(s)\,
ds-K^+(t)+K^-(t)+\int_0^t z(s)\, d W(s),\end{eqnarray*} where $
K^+\, (resp. K^-)$ is the weak (resp. strong) limit of
$\{K^{+,n}\}_{n=1}^\i \,(resp. \{K^{-,n}\}_{n=1}^\i)$ in $\cM^2$ and
$(K^+, K^-)\in \cN'^{\,2}\times \cN'^{\,2}$. Moreover, for any
$p\in[0,2)$,
$$ \lim_{n\ra \i} E\int_0^T |z^n(t)-z(t)|^p\, d t =0.$$ If
furthermore, $K^+$ is continuous, then we have $$ \lim_{n\ra \i} E\int_0^T
|z^n(t)-z(t)|^2 d t =0.$$ \el

 When applying the above generalized monotonic theorem, we need to compare
the increment of the increasing process $K_i^-$ for $i\in \L$. However, such
a kind of consideration does not seem to be available in the literature due
to the appearance of the lower barrier, which is implicit and thus varies
with the first unknown variable. The following lemma fills in such a gap,
which will be used in Section 3.

Assume that $\xi\in L^2(\O,\cF_T,P)$, and $L$ and $U$ are $\{\cF_t, 0\le t\le
T\}$-adapted continuous processes satisfying
$$ E[\sup_{0\leq t\leq T}\{|L(t)^+|^2+|U(t)^-|^2\}]<\i,\qq L(t)\leq U(t),\q  t\in [0,T].$$

Consider the following one-dimensional RBSDE with fixed double reflecting
barriers: \be\label{DRBSDE9}\left\{\ba{rcl}
\hat{Y}(t)&=&\displaystyle\xi+\int_t^T\psi(s,\hat{Y}(s),\hat{Z}(s))d
s-\int_t^T d
\hat{K}^-(s)\\[3mm]&&\displaystyle+\int_t^T d \hat{K}^+(s)-\int_t^T \hat{Z}(s)d W(s), \q t\in [0,T];\\[3mm]
L(t)&\le& \hat{Y}(t)\leq U(t),\q t\in [0,T];\\[3mm]
\displaystyle\int_0^T(\hat{Y}(s)&-&L(s))d\hat{K}^+(s)=0,\,\,
\displaystyle\int_0^T(\hat{Y}(s)-U(s))d \hat{K}^-(s)=0.\ea\right.\ee

\bde A barrier $S$ is called super-regular if there exists a sequence of
processes $\{S^n\}_{n=1}^{\i}$ such that

(i) $S^n(t)\geq S^{n+1}(t)$ and $\lim\limits_{n\ra\i}S^n(t)=S(t)$ for $t\in
[0, T]$;

(ii) For $n\geq 1$ and $t\in [0, T]$, we have $$d S^n(t)=u_n(t)\, dt+v_n(t)\,
d W(t)$$ where $u_n$ is an $\{\cF_t, 0\le t\le T\}$-adapted process such that
$$\sup\limits_{n\geq 1}\sup\limits_{0\leq t\leq T}|u_n(t)|<\i \q \hbox{\rm
and} \q v_n\in \cM^2.$$

A barrier $V$ is called sub-regular if the  barrier $-V$ is super-regular.
\ede

Note that the concept of our super-regular barrier is identical to the
definition of the regular upper barrier by  Hamad\`ene et al.
\cite{Hamadene1}.

\bl\label{l1} Assume that $\psi^1$ and $\psi^2$ satisfy Hypothesis 1 and that
the  barrier $U$ is super-regular. Assume that $\xi^1, \xi^2\in
L^2(\O,\cF_T,P)$, and $L^1, L^2$ and $U$ are $\{\cF_t, 0\le t\le T\}$-adapted
continuous processes satisfying
$$ E[\sup_{0\leq t\leq T}\{|L^j(t)^+|^2+|U(t)^-|^2\}]<\i,\qq L^j(t)\leq U(t),\q  t\in [0,T], \q j=1,2.$$
For $j=1,2,$ let $(\hat{Y}^j,\hat{Z}^j,\hat{K}^{+,j},\hat{K}^{-,j})$ be the
unique adapted solution of RBSDE~(\ref{DRBSDE9}) associated with data
$(\xi^j,\psi^j,L^j,U)$. Moreover, assume that

(i) $ \xi^1\leq \xi^2$;

(ii) $\psi^1(t,y,z)\leq \psi^2(t,y,z),\q \forall \,(y,z)\in R\times R^d$;

(iii) $L^1\leq L^2.$

Then we have

(1) $\hat{Y}^1(t)\leq \hat{Y}^2(t), \q 0\leq t\leq T $;

(2) $\hat{K}^{-,1}(r)-\hat{K}^{-,1}(s)\leq
\hat{K}^{-,2}(r)-\hat{K}^{-,2}(s),\q  0\leq s\leq r\leq T.$ \el

\begin{proof} For $j=1,2,$ and $n\geq 1$, the following penalized RBSDEs with a single
reflecting barrier: \begin{eqnarray*}\label{RBSDE1}\left\{\ba{rcl}
\hat{Y}^{j,n}(t)&=&\displaystyle\xi^j+\int_t^T\psi^j(s,\hat{Y}^{j,n}(s),\hat{Z}^{j,n}(s))\,
d
s -n\int_t^T(\hat{Y}^{j,n}(s)-U(s))^+\, ds\\[3mm]&&\displaystyle+\int_t^T
d\hat{K}^{j,n}(s)-\int_t^T \hat{Z}^{j,n}(s)\, d W(s),\q t\in [0,T];\\[3mm]
\hat{Y}^{j,n}(t)&\geq & L^j(t),\q t\in [0,T];\\[3mm]
\displaystyle\int_0^T(\hat{Y}^{j,n}(s)&-&L^j(s))\, d
\hat{K}^{j,n}(s)=0\ea\right.\end{eqnarray*} has a unique adapted
solution, denoted by $(\hat{Y}^{j,n}, \hat{Z}^{j,n},\hat{K}^{j,n})$.
In view of the comparison theorem \cite[Theorem 4.1]{El Karoui}, we
have
$$\hat{Y}^{1,n}(t)\leq \hat{Y}^{2,n}(t),\quad  0\leq t\leq T.$$
Noting that the barrier $U$ is super-regular, by the proof of \cite[Theorem
42.2 and Remark 42.3]{Hamadene1}, we have
$$\hat{Y}^j(t)=\lim_{n\ra\i}\hat{Y}^{j,n}(t),\q
\hat{K}^{-,j}(t)=\lim_{n\ra\i}n\int_0^t(\hat{Y}^{j,n}(s)-U(s))^+\, d s, \q
t\in [0,T].$$
 The desired results  then follow.\qquad\end{proof}

\br In a symmetric way, assuming that the lower barrier is sub-regular and
fixed, we can compare the increment of $\hat{K}^+$ when the upper barrier
varies.\er

The following lemma gives the continuous dependence of RBSDE with a r.c.l.l.
(right continuous with left limit) reflecting barrier, which will be used to
prove the lower minimal boundary condition.

\bl\label{l2} Assume that $\xi^j\in L^2(\O,\cF_T, P)$, $\psi^j$
satisfies Hypothesis 1, and $L^j\in\cS'^2$ for $j=1,2$.  For
$j=1,2$, denote by $(Y^j,
Z^j,K^j)\in\cS'^2\times\cM^2\times\cN'^{\,2}$ the unique adapted
solution of the following RBSDE: \begin{eqnarray*}\left\{\ba{rcl}
Y^j(t)&=&\displaystyle\xi^j+\int_t^T\psi^j(s,Y^j(s),Z^j(s))ds+\int_t^T
d
K^j(s)\\&&\displaystyle-\int_t^T Z^j(s)d W(s),\q t\in [0,T];\\[3mm] Y^j(t)&\geq& L^j(t),\q t\in [0,T];\\[3mm]\displaystyle \int_0^T
(Y^j(s-)&-&L^j(s-))d K^j(s)=0.\ea\right.\end{eqnarray*}  Set
$$\Delta
Y^j(t)\triangleq Y^j(t)-Y^j(t-),\q t\in [0,T]$$ for $j=1,2$. Then
there is a constant $c>0$ such that \begin{eqnarray*}\ba{rcl}&&
\displaystyle E\(\sup\limits_{0\leq t\leq
T}|Y^1(t)-Y^2(t)|^2+\int_0^T |Z^1(s)-Z^2(s)|^2d
s\)\\[3mm]&&\displaystyle+E\(|K^1(T)-K^2(T)|^2+\sum\limits_{0\leq t\leq
T}|\Delta Y^1(t)-\Delta Y^2(t)|^2\)\\[3mm]&\leq&\displaystyle c
E\(|\xi^1-\xi^2|^2+\int_0^T|\psi^1(s,Y^1(s),Z^1(s))-\psi^2(s,Y^1(s),Z^1(s))|^2
ds\)\\[3mm]&&+ c (E[\sup\limits_{0\leq t\leq
T}|L^1(t)-L^2(t)|^2])^{\frac{1}{2}}\Phi(T)^{\frac{1}{2}},\ea\end{eqnarray*}
where
$$\displaystyle\Phi(T)\triangleq\sum_{j=1}^2E\left[|\xi^j|^2+\int_0^T|\psi^j(t,0,0)|^2d
t+\sup_{0\leq t\leq T}|L^j(t)^+|^2\right].$$  \el

\begin{proof} The proof is similar to \cite[Proposition 3.6]{El Karoui} and is
omitted. \qquad\end{proof}

\br\label{r3} Lemma \ref{l2} can be extended to the multi-dimensional case
where for $ j=1,2,\, Y^j, \xi^j, \psi^j,K^j, S^j$ and $L^j$ are all
$R^m$-valued, $Z^j$ is $R^{m\times d}$-valued, and $ |z|\triangleq \sqrt{
trace(zz^T)}$ for $z\in R^{m\times d}$.\er

\section{Existence: the method of Picard iteration}

We have the following existence result for RBSDE~(\ref{DRBSDE}).
 \bt\label{t1}
Let Hypotheses 1, 2 and  3 be satisfied. Assume that the upper barrier $S$ is
super-regular with $S(t)\in Q(t)$ for $t\in [0,T]$, and that the terminal
value $\xi\in L^2(\O,\cF_T,P;R^m)$ takes values in $Q(T)$. Then
RBSDE~(\ref{DRBSDE}) has an adapted solution $(Y,Z,K^+,K^-)\in
(\cS^2)^m\times (\cM^2)^{m\times d}\times (\cN^2)^{2m}$.\et

\begin{proof} We use the method of Picard iteration. The whole proof is divided into
the following six steps.

\ms {\bf Step 1. Construction of Picard sequence of solutions $\{Y_i^n,
Z_i^n, K_i^n; i\in \L\}_{n\geq 0}.$}

For $ i\in \L,$  the following
 RBSDE with a single reflecting barrier:
\be\label{RBSDE2}\left\{\ba{rcl}
Y_i^0(t)&=&\displaystyle\xi_i+\int_t^T\psi(s,Y_i^0(s),Z_i^0(s),i)\, d
s-\int_t^T d
K_i^0(s)\\[3mm]&&\displaystyle-\int_t^T Z_i^0(s)\, d W(s),\q t\in [0,T];\\[3mm]
Y_i^0(t)&\leq &S_i(t),\q t\in [0,T];\\[3mm]\displaystyle \int_0^T(Y_i^0(s)&-&S_i(s))\, d
K_i^0(s)=0\ea\right.\ee has a unique adapted solution, denoted by $(Y^0_i,
Z_i^0, K_i^0)$.

For $n\geq 1$ and $i\in\L$, consider the following RBSDE with double
barriers: \be\label{DRBSDE2}\left\{\ba{rcl}
Y_i^n(t)&=&\displaystyle\xi_i+\int_t^T\psi(s,Y_i^n(s),Z_i^n(s),i)\, d
s-\int_t^T d
K_i^{-,n}(s)\\[3mm]&&\displaystyle+\int_t^T d K_i^{+,n}(s)-\int_t^T Z_i^n(s)\, d W(s),\q t\in [0,T];\\[3mm]
S_i(t)&\ge& Y_i^n(t)\geq \max\limits_{j\neq i, j\in \L}h_{i,j}(t,Y_j^{n-1}(t)),\q t\in [0,T];\\[3mm]
\displaystyle\int_0^T(Y_i^n(s)&-&\max\limits_{j\neq i, j\in
\L}h_{i,j}(s,Y_j^{n-1}(s)))\, dK_i^{+,n}(s)=0,\\[3mm]
\displaystyle\int_0^T(Y_i^n(s)&-&S_i(s))\, d K_i^{-,n}(s)=0.\ea\right.\ee
 Note that
$$
\max\limits_{j\neq i, j\in \L}h_{i,j}(t,Y_j^{n-1}(t))\le \max\limits_{j\neq
i, j\in \L}h_{i,j}(t,S_j(t))\le S_i(t), \q t\in [0,T],
$$
due to Hypothesis 2 and the assumption that $S(t)\in Q(t)$ for $t\in [0,T]$.
In view of~\cite[Theorem 42.2 and Remark 42.3]{Hamadene1},
RBSDE~(\ref{DRBSDE2}) has a unique solution
$$(Y_i^n,Z_i^n,K_i^{+,n},K_i^{-,n})\in
\cS^2\times\cM^2\times\cN^{\,2}\times\cN^{\,2}.$$

\ms {\bf Step 2. Convergence of $\{Y_i^n, K_i^{-,n}; i\in\L\}_{n\geq 1}.$}

RBSDE~(\ref{RBSDE2}) can be viewed as having the lower barrier $-\infty$.
Then from Lemma \ref{l1}, we have
$$ Y_i^0(t)\leq
Y_i^1(t), \q 0\leq t\leq T.$$
 Since $h_{i,j}(t,y)$ is increasing
in $y$, using Lemma \ref{l1} again, we know that the lower barrier
of RBSDE~(\ref{DRBSDE2}) is increasing with $n$ by induction.
Therefore, for $n\geq 1$ and $ i\in\L$, \begin{eqnarray*}\ba{l}
Y_i^n(t)\leq Y_i^{n+1}(t), \q K_i^{-,n}(t)\leq
K_i^{-,n+1}(t), \q 0\leq t\leq T;\\[3mm]K_i^{-,n}(r)-K_i^{-,n}(s)\leq
K_i^{-,n+1}(r)-K_i^{-,n+1}(s), \q 0\leq s\leq r\leq
T.\ea\end{eqnarray*} Hence, the sequence $ \{(Y_i^n(t),
K_i^{-,n}(t))\}_{n=1}^{\i}$ has a limit, denoted by $(Y_i(t),
K_i^-(t))$. From the last inequality, we have
\begin{eqnarray}\label{referee}
   K_i^{-,n}(r)-K_i^{-,n}(s)&\leq&
K_i^{-}(r)-K_i^{-}(s), \q 0\leq s\leq r\leq T,
\end{eqnarray}
which means that the process $K_i^-(s)-K_i^{-,n}(s), s\in [0,T]$ is
increasing. Since
$$ Y_i^0(t)\leq Y_i^n(t)\leq S_i(t), \q \forall \,t\in [0,T],\q i\in\L,$$
in view of the dominated convergence theorem, we have \be\label{8}
\lim_{n\ra\i}E\int_0^T|Y_i^n(t)-Y_i(t)|^2dt=0,\q i\in\L,\ee \be\label{3}
\sup_{n\ge 1}E[\sup_{0\leq t\leq T}|Y_i^n(t)|^2]\leq E\sup_{0\leq t\leq
T}(|Y_i^0(t)|^2+|S(t)|^2)<\i,\q i\in\L,\ee and \be Y_i^0(t)\leq Y_i(t)\leq
S_i(t), \q \forall \,t\in [0,T],\q i\in\L\ee which implies \be E[\sup_{0\leq
t\leq T}|Y_i(t)|^2]\leq E\sup_{0\leq t\leq T}(|Y_i^0(t)|^2+|S(t)|^2)<\i,\q
i\in\L.\ee

For $i\in\L$,  the following RBSDE \be\label{DRBSDE3}\left\{\ba{rcl}
\tilde{Y}_i(t)&=&\displaystyle\xi_i+\int_t^T\psi(s,\tilde{Y}_i(s),\tilde{Z}_i(s),i)\,
d s-\int_t^T d
\tilde{K}_i^-(s)\\[3mm]&&\displaystyle+\int_t^T d \tilde{K}_i^+(s)-\int_t^T \tilde{Z}_i(s)\, d W(s),\q t\in [0,T];\\[3mm]
S_i(t)&\ge&\tilde{Y}_i(t)\geq \max\limits_{j\neq i, j\in \L}h_{i,j}(t,S_j(t)),\q t\in [0,T];\\[3mm]
\displaystyle\int_0^T(\tilde{Y}_i(s)&-&\max\limits_{j\neq i, j\in
\L}h_{i,j}(s,S(s)))\, d\tilde{K}_i^+(s)=0,\\[3mm]
\displaystyle\int_0^T(\tilde{Y}_i(s)&-&S_i(s))\, d
\tilde{K}_i^-(s)=0\ea\right.\ee has a unique adapted solution, denoted by
$(\tilde{Y}_i,\tilde{Z}_i,
\tilde{K}_i^+,\tilde{K}_i^-)\in\cS^2\times\cM^2\times\cN^{\,2}\times\cN^{\,2}$.
 By Lemma \ref{l1}, we know that
for $n\geq 1$ and $i\in\L$, \begin{eqnarray*}\ba{l} K_i^{-,n}(t)\leq
K_i^{-,n+1}(t)\leq\tilde{K}_i^-(t), \q 0\leq t\leq T,\\[3mm]
K_i^{-,n}(r)-K_i^{-,n}(s)\leq K_i^{-,n+1}(r)-K_i^{-,n+1}(s)\leq
\tilde{K}_i^-(r)-\tilde{K}_i^-(s), \q 0\leq s\leq r\leq
T.\ea\end{eqnarray*} Then, $K_i^-$ is continuous for $i\in\L$.
Hence, \be\label{4} \qq \sup_{n\geq 1}E [\sup_{0\leq t\leq
T}|K_i^{-,n}(t)|^2]\leq  E [\sup_{0\leq t\leq T}|K_i^-(t)|^2]\leq E
[\sup_{0\leq t\leq T}|\tilde{K}_i^-(t)|^2]<\i,\,i\in\L.\ee  From
Dini's theorem, we have $$ \lim_{n\ra\i}\sup_{0\leq t\leq
T}|K_i^{-,n}(t)-K_i^-(t)|^2=0,\q i\in\L. $$ From the dominated
convergence theorem, we have \be\label{9}\lim_{n\ra\i}E[\sup_{0\leq
t\leq T}|K_i^{-,n}(t)-K_i^-(t)|^2]=0,\q i\in\L.\ee

\ms {\bf Step 3. Uniform boundedness of $\{\psi(\cdot, Y_i^n, Z_i^n,i),
Z_i^n, K_i^{+,n}; i\in \L\}_{n\geq 1}$ in $\cM^2\times \cM^2\times
\cN^{\,2}$.}

Applying It\^{o}'s lemma to $|Y_i^n(t)|^2$, we have for $i\in\L$,
\begin{eqnarray*}\ba{rcl}
&&\displaystyle|Y_i^n(t)|^2+\int_t^T |Z_i^n(s)|^2d s\\[3mm]
&=&\displaystyle\xi_i^2+2\int_t^T Y_i^n(s)\psi(s, Y_i^n(s),Z_i^n(s),i)d
s+2\int_t^T Y_i^n(s)d
K_i^{+,n}(s)\\[3mm]&&\displaystyle-2\int_t^T Y_i^n(s)d
K_i^{-,n}(s)-2\int_t^T Y_i^n(s)Z_i^n(s) d W(s).\ea\end{eqnarray*}
Using the Lipschitz property of $\psi$, the upper minimal boundary
condition in (\ref{DRBSDE2}) and the elementary inequality:
$\displaystyle ab\leq \frac{1}{\alpha}a^2+\alpha b^2,$ we have for
any arbitrary positive real number $\alpha$, \be\label{1}\ba{rcl}
&&\displaystyle E|Y_i^n(t)|^2+E[\int_t^T |Z_i^n(s)|^2d
s]\\[3mm]&\leq &\displaystyle E\(\xi_i^2 + 2\int_t^T Y_i^n(s)(\psi(s,0,0,i)
+C|Y_i^n(s)|+C|Z_i^n(s)|)d s\)\\[3mm]&&\displaystyle+E \(2\int_t^T Y_i^n(s)d
K_i^{+,n}(s) -2\int_t^T S_i(s)d K_i^{-,n}(s)\)\\[3mm]&\leq &\displaystyle
E\(\xi_i^2+\int_t^T |\psi(s,0,0,i)|^2d s+c\int_t^T |Y_i^n(s)|^2 d
s+\frac{1}{3}\int_t^T |Z_i^n(s)|^2d s\)\\[3mm]&&\displaystyle+E\(\frac{1}{\alpha}
[\sup_{0\leq t\leq T}|Y_i^n(t)|^2]+\alpha |K_i^{+,n}(T)|^2+\sup_{0\leq t\leq
T}|S_i(t)|^2+ |K_i^{-,n}(T)|^2\).\ea\ee Here and in the sequel, $c$ is a
positive constant whose value only depends on the Lipschitz coefficient $C$
and may change from line to line.

From RBSDE~(\ref{DRBSDE2}), we know that for
$i\in\L$,\begin{eqnarray*}
K_i^{+,n}(T)=Y_i^n(0)-\xi_i-\int_0^T\psi(s,Y_i^n(s),Z_i^n(s),i)d
s+K_i^{-,n}(T)+\int_0^T Z_i^n(s)d W(s).\end{eqnarray*} Hence, \be
\label{2} \qq \q E|K_i^{+,n}(T)|^2\leq c \(1+E\int_0^T
(|Y_i^n(s)|^2+|Z_i^n(s)|^2)ds+ E|K_i^{-,n}(T)|^2\),\q i\in\L. \ee
Substituting (\ref{2}) into (\ref{1}) and letting
$\alpha=\frac{1}{3c}$ and $t=0$, we have
$$E\int_0^T|Z_i^n(s)|^2d s\leq c
E\(1+\sup_{0\leq t\leq T}|Y_i^n(t)|^2+\int_0^T|Y_i^n(s)|^2d
s+|K_i^{-,n}(T)|^2\),\,i\in\L.$$ From (\ref{3}) and (\ref{4}), we know
\be\label{5.8}\sup_{n\geq 1}E\int_0^T|Z_i^n(s)|^2d s<\i,\q i\in\L.\ee Then
from (\ref{2}), we know \be \sup_{n\geq 1}E|K_i^{+,n}(T)|^2<\i,\q i\in\L.\ee
From (\ref{3}), (\ref{5.8}) and the Lipschitz property of $\psi$, we know
$$\sup_{n\geq 1} E\int_0^T |\psi(s, Y_i^n(s),Z_i^n(s),i)|^2 d s
<\i,\q i\in\L.$$ Therefore, without loss of generality, we can assume that
for $i\in\L$, $\{\psi(\cdot,Y_i^n,Z_i^n,i)\}_{n\geq 0}, \{Z_i^n\}_{n\geq 0},$
and $\{K_i^{+,n}\}_{n\geq 0}$ converge weakly in $\cM^2$ to $\psi_i, Z_i, $
and $ K_i^+$
 respectively.

\ms {\bf Step 4. Verification of the first equation of RBSDE (\ref{DRBSDE}).}

From the first equation of (\ref{DRBSDE2}), we have
$$ Y_i^n(t)=Y_i^n(0)-\int_0^t \psi(s, Y_i^n(s),Z_i^n(s),i)\, d
s-K_i^{+,n}(t)+K_i^{-,n}(t)+\int_0^t Z_i^n(s)\, d W(s).$$ All the
assumptions of the generalized monotonic limit theorem (see
Lemma~\ref{T2.1}) are shown to be satisfied in previous steps.
Therefore, for $i\in\L$, the limit $Y_i$ is r.c.l.l. and has the
form: \begin{eqnarray*} Y_i(t)=\xi_i+\int_t^T \psi_i(s)\, d
s-\int_t^T dK_i^-(s)+\int_t^T d K_i^+(s)-\int_t^T Z_i(s)\, d
W(s),\end{eqnarray*} and $K_i^+\in\cN'^{\,2}$. Moreover, for any
$p\in [0,2)$,$$\lim_{n\ra\i}E\int_0^T |Z_i^n(s)-Z_i(s)|^p \, ds=0,\q
i\in\L.$$ Hence, we have for $i\in\L$,
$$\lim_{n\ra\i}E \int_0^T|\psi(s,Y_i^n(s),Z_i^n(s),i)-\psi(s, Y_i(s),
Z_i(s),i)|^p\,  d s =0;$$
$$\psi_i(s)=\psi(s,Y_i(s),Z_i(s),i),\q a.e., a.s.;$$ \be\label{11}\ba{rcl}
Y_i(t)&=&\displaystyle\xi_i+\int_t^T \psi(s,Y_i(s),Z_i(s),i)\, d
s-\int_t^T dK_i^-(s)\\[3mm]&&\displaystyle+\int_t^T d
K_i^+(s)-\int_t^T Z_i(s)\, d W(s).\ea\ee

\ms {\bf Step 5. The upper and lower minimal boundary conditions.}

In view of RBSDE~(\ref{DRBSDE2}), we have
 $$\max\limits_{j\neq i,
j\in\L}\{h_{i,j}(t,Y_j^{n-1}(t))\}\leq Y_i^n(t)\leq S_i(t),\q t\in [0, T],
i\in \L.$$ Passing to the limit, we have \be\label{12}\max\limits_{j\neq i,
j\in\L}\{h_{i,j}(t,Y_j(t))\}\leq Y_i(t)\leq S_i(t),\q t\in [0, T], i\in
\L.\ee Since
$$\int_0^T(Y_i^n(s)-S_i(s))\, d K_i^{-,n}(s)=0\q\mbox{and}\q
Y_i^n(s)=Y_i^n(s-)\leq Y_i(s-)\leq S_i(s),$$ we have
$$0=\int_0^T(Y_i^n(s)-S_i(s))\, d K_i^{-,n}(s)\leq \int_0^T(Y_i(s-)-S_i(s))\, d K_i^{-,n}(s)\leq0,\q i\in\L.$$
Hence, $$\int_0^T(Y_i(s-)-S_i(s))\, d K_i^{-,n}(s)=0,\q i\in\L.$$ On the
other hand, in view of~(\ref{referee}), we have for $i\in\L$,
$$0\leq\int_0^T(S_i(s)-Y_i(s-))d (K_i^-(s)-K_i^{-,n}(s))\leq
\sup_{0\leq s\leq T}(S_i(s)-Y_i(s-))[K_i^-(T)-K_i^{-,n}(T)].$$ Since
$$\lim\limits_{n\ra \i} K_i^{-,n}(T)=K_i^-(T),\,i\in\L,$$ we have
$$\int_0^T (Y_i(s-)-S_i(s))\, d K_i^-(s)=\lim_{n\ra\i}\int_0^T (Y_i(s-)-S_i(s))\, d K_i^{-,n}(s)=0,\q
i\in\L.$$

We have just proved the upper minimal boundary condition. It remains to prove
the lower minimal boundary condition.  The technique used in \cite{Hamadene4}
is found difficult to be directly applied to our case since the corresponding
argument on the smallest $\psi$-supermartingale is not true in the case of
double barriers. We shall view the RBSDEs with double barriers as RBSDEs with
single lower barrier by taking the increasing processes $K^{-,n}$ as given.
For $i\in\L$ and $n\geq 1$,  the following RBSDE
\be\label{RBSDE4}\left\{\ba{rcl}
\bar{Y}_i^n(t)&=&\displaystyle\xi_i+\int_t^T\psi(s,\bar{Y}_i^n(s),\bar{Z}_i^n(s),i)d
s-K_i^{-,n}(T)\\[3mm]&&\displaystyle+K_i^{-,n}(t)+\int_t^T d
\bar{K}_i^{+,n}(s)-\int_t^T \bar{Z}_i^n(s)d
W(s),\\[3mm]\bar{Y}_i^n(t)&\geq &\wt h_i(t)\triangleq\max\limits_{j\neq i, j\in
\L}h_{i,j}(t,Y_j(t)), \q t\in [0,T],\\[3mm]
\displaystyle\int_0^T(\bar{Y}_i^n(s-)&-&\wt h_i(s-))\,
d\bar{K}_i^{+,n}(s)=0\ea\right.\ee has a unique adapted solution, denoted by
$(\bar{Y}_i^n,\bar{Z}_i^n,\bar{K}_i^{+,n})$.

Define $$\bar{X}_i^n\triangleq\bar{Y}_i^n-K_i^{-,n},\q
\psi_n(s,y,z,i)\triangleq \psi(s,y+K_i^{-,n}(s),z,i),\q i\in\L.$$ Then
$(\bar{X}_i^n,\bar{Z}_i^n,\bar{K}_i^{+,n})$ satisfies the following RBSDE:
\be\label{RBSDE5}\left\{\ba{rcl}
\bar{X}_i^n(t)&=&\displaystyle(\xi_i-K_i^{-,n}(T))+\int_t^T\psi_n(s,\bar{X}_i^n(s),\bar{Z}_i^n(s),i)\,
d
s\\[3mm]&&\displaystyle+\int_t^T d \bar{K}_i^{+,n}(s)-\int_t^T \bar{Z}_i^n(s)\, d
W(s),\\[3mm]\bar{X}_i^n(t)&\geq &\wt h_i(t)-K_i^{-,n}(t),\q t\in [0,T],\\[3mm]
\displaystyle\int_0^T\(\bar{X}_i^n(s-)&-&\wt h_i(s-)+K_i^{-,n}(s)\)\,
d\bar{K}_i^{+,n}(s)=0.\ea\right.\ee For $i\in\L$, let $
(\bar{Y}_i,\bar{Z}_i,\bar{K}_i^+)$ be the solution of the following RBSDE:
\be\label{RBSDE6}\left\{\ba{rcl}
\bar{Y}_i(t)&=&\displaystyle\xi_i+\int_t^T\psi(s,\bar{Y}_i(s),\bar{Z}_i(s),i)d
s-K_i^-(T)+K_i^-(t)\\[3mm]&&\displaystyle+\int_t^T d \bar{K}_i^+(s)-\int_t^T
\bar{Z}_i(s)d W(s),\q t\in [0,T];\\[3mm]\bar{Y}_i(t)&\geq &\wt h_i(t),\q t\in [0,T];\\[3mm]
\displaystyle\int_0^T(\bar{Y}_i(s-)&-&\wt h_i(s-))\,
d\bar{K}_i^+(s)=0.\ea\right.\ee

Define$$\bar{X}_i\triangleq\bar{Y}_i-K_i^-,\q\psi_-(s,y,z,i)\triangleq
\psi(s,y+K_i^-(s),z,i),\q i\in\L.$$ Then $(\bar{X}_i,\bar{Z}_i,\bar{K}_i^+)$
satisfies the following RBSDE: \be\label{RBSDE7}\left\{\ba{rcl}
\bar{X}_i(t)&=&\displaystyle(\xi_i-K_i^-(T))+\int_t^T\psi_-(s,\bar{X}_i(s),\bar{Z}_i(s),i)\,
d
s\\[3mm]&&\displaystyle+\int_t^T d \bar{K}_i^+(s)-\int_t^T \bar{Z}_i(s)\, d
W(s),\q t\in [0,T];\\[3mm]\bar{X}_i(t)&\geq &\wt h_i(t)-K_i^-(t),\q t\in [0,T];\\[3mm]
\displaystyle\int_0^T\(\bar{X}_i(s-)&-&\wt
h_i(s-)+K_i^-(s)\)d\bar{K}_i^+(s)=0.\ea\right.\ee
 Since
\begin{eqnarray*}\ba{rcl}&&\psi_n(s,y,z,i)-\psi_-(s,y,z,i)\\[3mm]&=&\psi(s,y+K_i^{-,n}(s),z,i)-\psi(s,y+K_i^-(s),z,i)\\[3mm]&\leq
& C|K_i^{-,n}(s)-K_i^-(s)|,\ea\end{eqnarray*}  in view of Lemma
\ref{l2}, we have \be\label{31}\ba{rcl} &&E[\sup\limits_{0\leq t\leq
T}|\bar{X}_i^n(t)-\bar{X}_i(t)|^2]\\&\leq &\displaystyle c
E[|K_i^-(T)-K_i^{-, n}(T)|^2+C\int_0^T |K_i^-(t)-K_i^{-, n}(t)|^2
dt]\\[3mm]&&\ds+c \{E[\sup_{0\leq t\leq
T}|K_i^-(t)-K_i^{-,
n}(t)|^2]\}^{\frac{1}{2}}(\Phi^n(T))^{\frac{1}{2}}, \ea\ee where
\begin{eqnarray*}\ba{rcl} \Phi^n(T)&\triangleq&\displaystyle E[ (\xi_i-K_i^{-,
n}(T))^2+\int_0^T|\psi_{n}(s, 0, 0, i)|^2d
s]\\[3mm]&&+E[\sup\limits_{0\leq t\leq T}((\wt h_i(t)-K_i^{-, n}(t))^+)^2+\sup\limits_{0\leq t\leq T}((\wt h_i(t)-K_i^-(t))^+)^2]\\[3mm]&&\displaystyle+E[(\xi_i-K_i^-(T))^2+\int_0^T|\psi_{-}(s, 0, 0, i)|^2d
s]. \ea\end{eqnarray*} Since
\begin{eqnarray*}\ba{l}\psi_n(s,0,0,i)\leq\psi(s,0,0,i)+C
K_i^{-,n}(s),\\[3mm]\psi_-(s,0,0,i)\leq\psi(s,0,0,i)+C K_i^-(s),\ea\end{eqnarray*}
\begin{eqnarray*}\ba{rcl}\sup\limits_{0\leq t\leq T}((\wt
h_i(t)-K_i^{-,n}(t))^+)^2&\leq& \sup\limits_{0\leq t\leq
T}((\wt h_i(t))^+)^2\\[3mm]&\leq& \sup\limits_{0\leq t\leq
T}((\max\limits_{j\neq i, j\in \L}Y_j(t))^+)^2
\\[3mm]&\leq &\displaystyle\sum\limits_{j\in\L}\sup\limits_{0\leq t\leq T}|Y_j(t)|^2,\ea\end{eqnarray*}
 and  $$\sup\limits_{0\leq t\leq T}((\wt h_i(t)-K_i^-(t))^+)^2\leq
\sum\limits_{j\in\L}\sup\limits_{0\leq t\leq T}|Y_j(t)|^2,$$ we have
\begin{eqnarray*}\ba{rcl} \Phi^n(T)&\leq&\displaystyle
4E(\xi_i^2+\int_0^T|\psi(s,0,0,i)|^2d
s)+2\sum\limits_{j\in\L}E[\sup\limits_{0\leq t\leq
T}|Y_j(t)|^2]\\[3mm]&&+2(C^2T+1)E(|K_i^{-,n}(T)|^2+|K_i^-(T)|^2).\ea\end{eqnarray*}
From (\ref{3}) and (\ref{4}), we see \be\label{38} \sup_{n\geq
1}\Phi^n(T)<\i.\ee From (\ref{9}), (\ref{31}) and (\ref{38}), we see
$$\lim_{n\ra\i}E[\sup_{0\leq t\leq T}|\bar{X}_i^n(t)-\bar{X}_i(t)|^2]=0,\q i\in\L.$$   So
there is a subsequence of $\{\bar{X}_i^n\}_{n\geq 1}$ converging to
$\bar{X}_i, a.e., a.s.$ Without loss of generality, assume that
\be\label{30}\lim_{n\ra\i}\bar{X}_i^n= \bar{X}_i,\q a.e., a.s.,\q i\in\L.\ee
Set
$$X_i^n\triangleq Y_i^n-K_i^{-,n},\q i\in\L.$$ Then
from reflected BSDEs (\ref{DRBSDE2}) we know that $(X_i^n, Z_i^n, K_i^{+,n})$
is the solution of the following reflected BSDEs with single reflecting
barrier: \be\label{RBSDE10}\left\{\ba{rcl}
X_i^n(t)&=&\displaystyle(\xi_i-K_i^{-,n}(T))+\int_t^T\psi_{n}(s, X_i^n(s),
Z_i^n(s),i)d s\\[3mm]&&\displaystyle+\int_t^T d K_i^{+,n}(s)-\int_t^T  Z_i^n(s)d
W(s),\q t\in [0,T];\\[3mm] X_i^n(t)&\geq &\max\limits_{j\neq i, j\in
\L}h_{i,j}(t,Y_j^{n-1}(t))-K_i^{-,n}(t), \q t\in [0,T];\\[3mm]
\displaystyle\int_0^T\(X_i^n(s)&-&\max\limits_{j\neq i, j\in
\L}h_{i,j}(s,Y_j^{n-1}(s))+K_i^{-,n}(s)\)d K_i^{+,n}(s)=0.\ea\right.\ee
Comparing it with reflected BSDEs (\ref{RBSDE5}) and using the comparison
theorem for r.c.l.l. reflecting barrier \cite[Theorem 1.5]{Hamadene2}, we
know that
$$ \bar{X}_i^n(t)\geq X_i^n(t)=Y_i^n(t)-K_i^{-,n}(t),\q  (t,i) \in [0,T]\times \L.$$ In view of (\ref{30}), we have \be\label{13}
\bar{X}_i(t)\geq Y_i(t)-K_i^-(t),\q  (t,i) \in [0,T]\times \L.\ee Note that
due to the appearance of the additional fixed upper barrier, it is not clear
whether the lower barrier of (\ref{RBSDE7}) is not less than that of
(\ref{RBSDE10}). Such a difficulty is got around by comparing (\ref{RBSDE10})
and (\ref{RBSDE5}).

On the other hand, from (\ref{RBSDE7}) and \cite[Theorem 2.1]{Lepeltier}, we
know that $\bar{X}_i(\cdot)$ is the smallest $\psi_-$-supermartingale with
the lower barrier $\{\wt h_i(t)-K_i^-(t),0\leq t\leq T\}.$ From (\ref{11})
and (\ref{12}), it can be easily obtained that $\{Y_i(t)-K_i^-(t),\,0\leq
t\leq T\}$ is a $\psi_-$-supermartingale with the same lower barrier and
terminal value. Hence, $$ \bar{X}_i(t)\leq Y_i(t)-K_i^-(t),\q \forall \,t\leq
T, i\in\L .$$ Together with (\ref{13}), we have
$$ \bar{X}_i(t)= Y_i(t)-K_i^-(t),\q  \forall \,t\leq T, i\in\L .$$  Then
$$\bar{Y}_i(t)=Y_i(t),\q  \forall \,t\leq T, i\in\L .$$
From the uniqueness of the Doob-Meyer Decomposition, it follows that
$$ \bar{Z}_i(t)=Z_i(t), \q \bar{K}_i^+(t)=K_i^+(t),\q \forall \,0\leq
t\leq T, i\in\L .$$ Then we have
\be\label{3.118}\int_0^T\left(Y_i(s-)-\max\limits_{j\neq i,\, j\in \L}h_{i,
j}(s, Y_j(s-))\right)dK_i^+(s)=0,\q\forall \, i\in\L.\ee Hence, for
$i\in\L,\q (Y_i, Z_i ,K_i^+, K_i^-)$ almost satisfies RBSDE~(\ref{DRBSDE})
except that both minimal boundary conditions are replaced by
$$\int_0^T\(Y_i(s-)-\max\limits_{j\neq i, j\in
\L}h_{i,j}(s,Y_j(s-))\)dK_i^+(s)=0,\q \int_0^T(Y_i(s-)-S(s))d K_i^-(s)=0.$$
If we further prove the continuity of $Y_i$ , then we know $\{(Y_i, Z_i
,K_i^+, K_i^-), i\in \L\}$ is an adapted solution of RBSDE~(\ref{DRBSDE}).

\ms {\bf Step 6. The time continuity of $Y$ and $K^+$.}

Since $K_i^-$ is continuous, we have $\Delta Y_i(t)=-\Delta
K_i^+(t)\leq 0.$ Since the fact that $\Delta Y_{i_1}(t^*)<0$ for
some $(i_1, t^*)\in\L\times [0,T]$,  contradicts Hypothesis 3, as
shown in~\cite[Theorem 2.4]{Hamadene4}, we have \be \Delta Y_i(t)=0,
\q 0\leq t\leq T, \,i\in\L,\ee that is, $Y$ is time continuous.
Hence $K_i^+$ is time continuous. \qquad\end{proof}

\section{Existence: the penalty method}
In Theorem \ref{t1}, every component of  the upper barrier S is assumed to be
super-regular. Let $k(i,j)$ be the switching cost from state i to state j in
the optimal switching problem (see Hu and Tang \cite{Hu2}) satisfying
Hypothesis 3', and let the function $h_{i,j}$ introduced in the preceding
section take the particular form $h_{i,j}(t,y)=y-k(i,j)$. Then we can prove
by a penalty method the existence of an adapted solution to
RBSDE~(\ref{DRBSDE}) without the super-regularity assumption on the upper
barrier $S$. Moreover, we shall avoid using Lemma~\ref{T2.1} since its proof
itself is very lengthy and further invokes very technical lemmas of
Peng~\cite{Peng0}.

\subsection{Multi-dimensional RBSDEs with  fixed single
reflecting barrier}

In what follows, we consider the multi-dimensional RBSDEs with  fixed single
reflecting barrier, show the existence and uniqueness by a penalty method,
and give a comparison theorem.

For two $m$-dimensional vectors $x\triangleq (x_1, \cdots, x_m)^T$ and
$y\triangleq (y_1, \cdots, y_m)^T$, by $x\le y$ we mean  $x_i\le y_i$ for
$i\in \L$. For a vector $x\triangleq (x_1, \cdots, x_m)^T$, $x^+$ is defined
as the $m$-dimensional vector $(x_1^+,\cdots,x_m^+)^T$.

Consider the following multi-dimensional RBSDE with  fixed single reflecting
barrier:
 \be\label{RBSDE8}\left\{\ba{rcl}
Y(t)&=&\displaystyle\xi+\int_t^T\phi(s,Y(s),Z(s))\, d s-\int_t^T d
K(s)\\[3mm]&&\displaystyle-\int_t^T Z(s)\, d W(s),\q t\in [0,T];\\[3mm]
Y(t)&\leq & S(t), \q t\in [0,T];\\[3mm]
\displaystyle\int_0^T(Y_i(s)&-&S_i(s))\, d K_i(s)=0, \q  i\in
\L.\ea\right.\ee We make the following assumption on the generator $\phi$,
the terminal value $\xi\triangleq (\xi_1,\cdots, \xi_m)^T$, and the barrier
$S\triangleq (S_1,\cdots, S_m)^T$.

\ms \textbf{Hypothesis 4.} (i) The process $\phi(\cdot,0,0)\in (\cM^2)^m$.
For $i\in \L$, $\xi_i\in L^2(\O,\cF_T,P)$ and $S_i\in \cS^2$ with $\xi_i\le
S_i(T).$

(ii) There is a constant $C>0$ such that for any $(t,y,y',z,z')\in
[0,T]\times (R^m)^2\times (R^{m\times d})^2,$ we have
\be\label{50}\ba{rcl}|\phi(t,y,z)-\phi(t,y',z')|&\leq&
C(|y-y'|+|z-z'|),\\[3mm] -4\langle
y^-,\phi(t,y^++y',z)-\phi(t,y',z')\rangle&\leq&\displaystyle
2\sum_{i=1}^m\chi_{\{y_i<0\}}|z_i-z'_i|^2+C|y^-|^2.\ea\ee
\ms We
have \bt\label{t4} Let Hypothesis 4 be satisfied. Then
RBSDE~(\ref{RBSDE8}) has a unique adapted solution $(Y,Z,K)\in
(\cS^2)^m\times (\cM^2)^{m\times d}\times (\cN^{\,2})^m$. \et

\begin{proof} For any  positive integer n, consider the following penalized BSDE:
\be\label{BSDE8}\ba{rcl} Y^n(t)&=&\displaystyle
\xi+\int_t^T\phi(s,Y^n(s),Z^n(s))d
s-n\int_t^T(Y^n(s)-S(s))^+\, d s\\[3mm]&&
\displaystyle-\int_t^T Z^n(s)\, d W(s), \q t\in [0,T]. \ea\ee Note that
$(Y^n(s)-S(s))^+$ has been defined component-wisely in the second paragraph
of this subsection. From Pardoux and Peng~\cite{Pardoux}, we know that BSDE
(\ref{BSDE8}) has a unique adapted solution $(Y^n, Z^n)\in
(\cS^2)^m\times(\cM^2)^{m\times d}$ for each $n$.

Define for $(t,y,z)\in [0,T]\times R^m\times R^{m\times d}$,
$$\phi^n(t,y,z)\triangleq \phi(t,y,z)-n(y-S(t))^+$$
and $$ K^n(t)\triangleq n\int_0^t(Y^n(s)-S(s))^+\, d s.$$

In view of Hypothesis 4 (ii), we have for all $y\in R^m$,
\begin{eqnarray*}\ba{rcl} && -4\langle y^-,
\phi^n(s,y^++y',z)-\phi^{n+1}(s,y',z')\rangle\\[3mm]&\leq&-4\langle
y^-,\phi(s,y^++y',z)-\phi(s,y',z')+(y'-S(s))^+\rangle\\[3mm]&\leq&-4\langle
y^-,\phi(s,y^++y',z)-\phi(s,y',z')\rangle\\[3mm]&\displaystyle\leq&\displaystyle 2\sum_{i=1}^m\chi_{\{y_i<0\}}|z_i-z'_i|^2+C|y^-|^2.\ea\end{eqnarray*}
Applying the comparison theorem of multi-dimensional BSDEs (see
\cite[Theorem 2.1]{Hu1}), we deduce that $$Y^{n+1}(t)\leq Y^n(t), \q
0\leq t\leq T.$$ For $t\in [0,T]$, the sequence
$\{Y^n(t)\}_{n\geq1}$ almost surely admits a limit, which is denoted
by $Y(t)$ below.

Applying It\^{o}'s lemma to compute $|Y^n(t)|^2$,  in view of the following
inequality
$$\int_t^T \langle Y^n(s), d K^n(s)\rangle \geq \int_t^T \langle S(s), d
K^n(s)\rangle,$$ it is more or less standard to derive the following (see
Gegout-Petit and Pardoux~\cite{Gegout-Petit}):
 \be\label{5}E|K^n(T)|^2+\sup_{0\leq t\leq T}E|Y^n(t)|^2 +E\int_0^T|Z^n(s)|^2ds\leq c, \q
n=1,2,\cdots\ee In view of (\ref{BSDE8}),  applying the
Burkholder-Davis-Gundy (abbreviated as BDG below) inequality, we have
$$E\sup_{0\leq t\leq T}|Y^n(t)|^2+E\int_0^T|Z^n(s)|^2\, ds+E|K^n(T)|^2\leq
c, \q n=1,2,\cdots$$ Recalling that $ Y(t)=\lim\limits_{n\ra\infty}Y^n(t),$
using Fatou's lemma, we have
$$E\sup_{0\leq t\leq T}|Y(t)|^2\leq  E\mathop{\underline{\lim}}_{n\to \infty}\sup_{0\leq t\leq T}|Y^n(t)|^2\leq  \mathop{\underline{\lim}}_{n\to \infty}E\sup_{0\leq t\leq T}|Y^n(t)|^2\leq c.$$
Since $|Y^n(s)-Y(s)|^2\leq |Y^0(s)-Y(s)|^2$ for $s\in [0,T]$ (noting that
$Y^n$ converges to $Y$ in a decreasing manner),  we have from Lebesgue's
dominated convergence theorem that
$$\lim_{n\ra\infty}E\int_0^T |Y^n(s)-Y(s)|^2d s=0.$$

For positive integers $n_1$ and $n_2$, applying It\^{o}'s lemma to
compute $|Y^{n_1}(t)-Y^{n_2}(t)|^2$, it is standard to derive the
following (see Gegout-Petit and Pardoux~\cite{Gegout-Petit}):
\begin{eqnarray*}\ba{rcl}&&\displaystyle
E|Y^{n_1}(0)-Y^{n_2}(0)|^2+E\int_0^T|Z^{n_1}(s)-Z^{n_2}(s)|^2d
s\\[3mm]&\leq&\displaystyle 2C(C+1) E\int_0^T|Y^{n_1}(s)-Y^{n_2}(s)|^2d
s+\frac{1}{2}E\int_0^T|Z^{n_1}(s)-Z^{n_2}(s)|^2d s\\[3mm]&&\displaystyle +2
E\int_0^T\langle (Y^{n_2}(s)-S(s))^+, d K^{n_1}(s)\rangle +2
E\int_0^T\langle (Y^{n_1}(s)-S(s))^+, d
K^{n_2}(s)\rangle.\ea\end{eqnarray*} As a consequence,
\be\label{81}\ba{rcl}\displaystyle
E\int_0^T|Z^{n_1}(s)-Z^{n_2}(s)|^2d s&\leq&\displaystyle 4C(C+1)
E\int_0^T|Y^{n_1}(s)-Y^{n_2}(s)|^2d s\\[3mm]&&\displaystyle+4
E\int_0^T\langle (Y^{n_2}(s)-S(s))^+, d K^{n_1}(s)\rangle \\[3mm]&&\displaystyle +4
E\int_0^T\langle (Y^{n_1}(s)-S(s))^+, d K^{n_2}(s)\rangle.\ea\ee

Proceed identically as in the proof of \cite[Lemma 6.1]{El Karoui}, we have
\be\label{88}E\sup_{0\leq t\leq T}|(Y^n(t)-S(t))^+|^2\ra 0 \mbox{\q as\q}
n\ra\infty.\ee

Then from (\ref{5}) and (\ref{88}), we know that as $n_1,n_2\ra\i$,
$$E\int_0^T\langle (Y^{n_2}(s)-S(s))^+, d K^{n_1}(s)\rangle +
E\int_0^T\langle (Y^{n_1}(s)-S(s))^+, d K^{n_2}(s)\rangle \ra 0.$$ Together
with (\ref{81}), we obtain
$$\lim_{n_1,n_2\ra\i} E\int_0^T|Z^{n_1}(s)-Z^{n_2}(s)|^2d s =0.$$
In view of~(\ref{BSDE8}), using the BDG inequality, we conclude
\begin{eqnarray*} \lim_{n_1,n_2\ra \infty}E\sup_{0\leq t\leq
T}|Y^{n_1}(t)-Y^{n_2}(t)|^2=0.\end{eqnarray*} In view of
(\ref{BSDE8}) and the definition of $K^n$, we know
$$\lim_{n_1,n_2\ra\infty}E\sup_{0\leq t\leq T}|K^{n_1}(t)-K^{n_2}(t)|^2=0.$$
From the above convergence,  we conclude that there exists  $( Z, K)\in
(\cM^2)^{m\times d}\times(\cN^{\,2})^m$ such that
$$\lim_{n\ra\infty}E\int_0^T|Z^n(s)-Z(s)|^2\, d
s=0 \hbox{ \rm and } \lim_{n\ra\infty}E\sup_{0\leq t\leq
T}|K^n(t)-K(t)|^2=0.$$ Passing to limit in equation~(\ref{BSDE8}), we know
that
$$(Y ,Z ,K) \in
(\cS^2)^m\times(\cM^2)^{m\times d}\times(\cN^{\,2})^m$$ satisfies the
following equation:
$$Y(t)=\xi+\int_t^T\phi(s,Y(s),Z(s))\, d s-\int_t^T d K(s)-\int_t^T
Z(s)\, d W(s), \q t\in [0,T].$$

Since $Y(t)\leq Y^n(t)$,   we have
$$E\sup_{0\leq t\leq T}|(Y(t)-S(t))^+|^2\leq E\sup_{0\leq t\leq T}|(Y^n(t)-S(t))^+|^2,\q n=1,2,\ldots,$$
and further from (\ref{88}),
$$E\sup_{0\leq t\leq T}|(Y(t)-S(t))^+|^2\leq \lim_{n\ra\infty}
E\sup_{0\leq t\leq T}|(Y^n(t)-S(t))^+|^2=0,$$ which implies that
\be\label{6}Y(t)\leq S(t),\q 0\leq t\leq T.\ee Therefore,
\be\label{6.1}\int_0^T(Y_i(t)-S_i(t))\, d K_i(t)\leq 0,\q i\in\L.\ee Since
$(Y^n, K^n)$ tends to $(Y,K)$ in $(\cS^2)^{2m}$,
 we have for $i\in\L$,
\be\label{90}\ba{rcl}\displaystyle\int_0^T(Y_i(t)-S_i(t))\, d
K_i(t)&=&\displaystyle\lim_{n\ra \i}\int_0^T(Y_i^n(t)-S_i(t))\, d
K_i^n(t)\\[3mm]&=&\displaystyle\lim_{n\ra
\i}n\int_0^T(Y_i^n(t)-S_i(t))(Y_i^n(t)-S_i(t))^+d t\\[3mm]&\geq&0.\ea\ee
Thus, \be\label{7}\int_0^T(Y_i(t)-S_i(t))\, d K_i(t)=0,\q i\in\L.\ee We
conclude that $(Y ,Z ,K)$ is an adapted  solution to RBSDE~(\ref{RBSDE8}).

Uniqueness of the solution follows from Lemma \ref{l2} and Remark
\ref{r3}.\qquad
\end{proof}

\br For the existence and uniqueness for multi-dimensional RBSDEs, we refer
the reader to Gegout-Petit and Pardoux~\cite{Gegout-Petit}. \er

\br The comparison theorem of multi-dimensional BSDEs is first established by
Hu and Peng \cite{Hu1} under the stronger conditions on the generator $\phi$
that  $\phi(\cdot,y,z)$ is continuous for any fixed $(y,z)$ and
$\phi(\cdot,0,0)\in (\cS^2)^m$. By the method of approximation, it can be
shown that the comparison theorem still holds if Hypothesis 4 (i) is
satisfied. \er

Thanks to the above existence and uniqueness result, we can prove the
following comparison theorem for multi-dimensional RBSDEs with
 fixed single reflecting barrier.

\bt\label{t7} Assume that $(\phi^1, \xi^1), (\phi^2, \xi^2),$ and $S$
satisfy Hypothesis 4. Further,  assume that

 (i) $\xi^1\geq \xi^2;$

 (ii) There is a positive constant $C'$ such that for $(y,y')\in (R^m)^2, (z,z')\in (R^{m\times d})^2,$ and $t\in[0,T]$,
\be\label{58}\q -4\langle
y^-,\phi^1(t,y^++y',z)-\phi^2(t,y',z')\rangle\leq
2\sum_{i=1}^m\chi_{\{y_i<0\}}|z_i-z'_i|^2+C'|y^-|^2.\ee For $j=1,2,$
denote by $(Y^j,Z^j,K^j)$ the adapted solution of
RBSDE~(\ref{RBSDE8}) associate with the data $(\xi^j, \phi^j, S)$.
Then, we have

(1) $Y^1(t)\geq Y^2(t), \q   0\leq t\leq T$;

(2) $K^1(r)-K^1(s)\geq K^2(r)-K^2(s), \q  0\leq s\leq r\leq T.$ \et

 \begin{proof} For $j=1,2$ and positive integer n,   the following BSDE: \be\label{BSDE9}\ba{rcl}
Y^{j,n}(t)&=&\displaystyle \xi^j+\int_t^T\phi^j(s,Y^{j,n}(s),Z^{j,n}(s))d
s-n\int_t^T(Y^{j,n}(s)-S^j(s))^+\, d s\\[3mm]&&
\displaystyle-\int_t^T Z^{j,n}(s)\, d W(s), \q t\in [0,T] \ea\ee has a unique
adapted solution, denoted by $(Y^{j,n},Z^{j,n},K^{j,n})$.
  In view of
(\ref{58}), we have \begin{eqnarray*}\ba{rcl}&&-4\langle
y^-,(\phi^1(t,y^++y',z)-n(y^++y'-S(t))^+)
-(\phi^2(t,y',z')-n(y'-S(t))^+)\rangle\\[3mm]&\leq&-4\langle
y^-,\phi^1(t,y^++y',z) -\phi^2(t,y',z')-n y^+\rangle \\[3mm]&\leq
&\displaystyle
2\sum_{i=1}^m\chi_{\{y_i<0\}}|z_i-z'_i|^2+C'|y^-|^2.\ea\end{eqnarray*}
By the comparison theorem of multi-dimensional BSDEs (see
\cite[Theorem 2.1]{Hu1}), it follows that
\be\label{59}Y^{1,n}(t)\geq Y^{2,n}(t), \q t\in[0,T],
n=1,2,\cdots.\ee Then from the proof of Theorem \ref{t4}, we know
that for $t\in[0,T]$ and $j=1,2,$ \be\label{60}\ba{rcl}
Y^j(t)&=&\displaystyle\lim_{n\ra\i}Y^{j,n}(t),\\[3mm]
K^j(t)&=&\displaystyle \lim_{n\ra\i}n\int_0^t(Y^{j,n}(s)-S(s))^+d
s.\ea\ee The desired results then follow from (\ref{59}) and
(\ref{60}). \qquad\end{proof}

{\it Note added after the second round of review}: We thank the reviewer in
the second round of review who reminds us of the recent work of Wu and
Xiao~\cite{WuXiao}. On one hand, Wu and Xiao~\cite[Theorem 2.2, page
1821]{WuXiao} is a more general existence and uniqueness result than our
Theorem~\ref{t4} since the second inequality in~(\ref{50}) of our Hypothesis
4 (ii) is not required in their theorem. On the other hand, concerning the
comparison theorem for multidimensional RBSDEs, in addition to comparison of
the first component of the solution,  our Theorem~\ref{t7} further includes
comparison on the third component (the so-called increasing process) of the
solution, which is missing in Wu and Xiao~\cite[Theorem 3.1, page
1825]{WuXiao}.

\subsection{Multi-dimensional RBSDEs with an oblique reflection}

Consider the following RBSDE: for $i\in \L$,
\be\label{DRBSDE8}\left\{\ba{rcl}
Y_i(t)&=&\displaystyle\xi_i+\int_t^T\psi(s,Y_i(s),Z_i(s),i)d s-\int_t^T d
K_i^-(s)\\[3mm]&&\displaystyle+\int_t^T d K_i^+(s)-\int_t^T Z_i(s)d W(s), \q t\in [0,T];\\[3mm]
S_i(t)&\ge&Y_i(t)\geq \max\limits_{j\neq i, j\in \L}\{Y_j(t)-k(i,j)\},\q t\in [0,T];\\[3mm]
&&\displaystyle\int_0^T\left(Y_i(s)-\max\limits_{j\neq i, j\in
\L}\{Y_j(s)-k(i,j)\}\right)dK_i^+(s)=0;\\[3mm]
&&\displaystyle\int_0^T(Y_i(s)-S_i(s))d K_i^-(s)=0.\ea\right.\ee

The following theorem presents the existence of the solution without
super-regularity assumption on $S$.

\bt\label{t6} Let Hypotheses 1 and 3' be satisfied.  Assume that $S\in
(\cS^2)^m$ with $S(t)\in \tilde Q(t)$ for $t\in [0,T]$ and $\xi\in
L^2(\O,\cF_T,P; R^m)$ with $\xi(\o)\in \tilde{Q}(T)$. Here we have defined
for $t\in [0,T],$
$$\tilde{Q}(t)\triangleq\{(y_1,\cdots, y_m)^T\in R^m :
y_j-k(i,j) \leq y_i\leq S_i(t) ,\forall \,i,j\in \L, j\neq i\}.$$ Then
RBSDE~(\ref{DRBSDE8}) has an adapted  solution $(Y,Z,K^+,K^-)\in
(\cS^2)^m\times (\cM^2)^{m\times d}\times (\cN^{\,2})^{2m}$.\et

\begin{proof}  The proof is divided into four steps.

\ms {\bf Step 1. The approximating sequence of penalized RBSDEs.}

For any positive integer $n$, consider the following  RBSDE: $\forall
\,i\in\L$,\be\label{RBSDE9}\left\{\ba{rcl}
Y_i^n(t)&=&\displaystyle\xi_i+\int_t^T\psi(s,Y_i^n(s),Z_i^n(s),i)\, d
s-\int_t^T d
K_i^{-,n}(s)\\[3mm]&&\displaystyle+n\sum_{l=1}^m\int_t^T(Y_i^n(s)-Y_l^n(s)+k(i,l))^-\, d s\\
&&\displaystyle -\int_t^T Z_i^n(s)\, d W(s),\q t\in [0,T];\\[3mm]
Y_i^n(t)&\leq &S_i(t),\q t\in [0,T];\\[3mm]
\displaystyle\int_0^T(Y_i^n(s)&-&S_i(s))d K_i^{-,n}(s)=0.\qq\qq\qq
\ea\right.\ee These RBSDEs satisfy Hypothesis 4.

In fact, define for $(t,y,z,i)\in [0,T]\times R^m\times R^{m\times
d}\times \L$,
\begin{eqnarray*}\ba{rcl}\bar{\psi}^n(t,y,z,i)&\triangleq
&\displaystyle\psi(t,y_i,z_i,i)+n\sum_{l=1}^m(y_i-y_l+k(i,l))^-,\\[3mm]\bar{\psi}^n(t,y,z)&\triangleq
&
(\bar{\psi}^n(t,y,z,1),\cdots,\bar{\psi}^n(t,y,z,m))^T.\ea\end{eqnarray*}
Since \be\ba{rcl} &&\langle
y_i^-,(y_i^++y'_i-y_l^+-y'_l+k(i,l))^--(y'_i-y'_l+k(i,l))^-\rangle\\[3mm]&=&
\langle y_i^-, (y'_i-y_l^+-y'_l+k(i,l))^--(y'_i-y'_l+k(i,l))^-\rangle\geq
0\ea\ee for  $y,y'\in R^m$ and $ i, l \in \L$, and $\bar{\psi}^n(t,y,z,i)$
does not depend  on $z_j$ for $j\not=i$, we have for $(y,y')\in (R^m)^2,
(z,z')\in (R^{m\times d})^2,$ and $(i, l) \in (\L)^2,$
\be\label{89}\ba{rcl}&&-4\langle
y^-,\bar{\psi}^n(t,y^++y',z)-\bar{\psi}^n(t,y',z')\rangle\\[3mm]&=&\displaystyle-4\sum_{i=1}^m\langle y_i^-,\psi(t,y_i^++y'_i,z_i,i)-\psi(t,y'_i,z'_i,i)\rangle\\[3mm]&&\displaystyle-4n\sum_{i,\,l=1}^m\langle y_i^-, (y_i^++y'_i-y_l^+-y'_l+k(i,l))^--(y'_i-y'_l+k(i,l))^-\rangle\\[3mm]&\leq &\displaystyle -4\sum_{i=1}^m\langle y_i^-,\psi(t,y_i^++y'_i,z_i,i)-\psi(t,y'_i,z'_i,i)\rangle\\[3mm]&\leq&
\displaystyle2\sum_{i=1}^m\chi_{\{y_i<0\}}|z_i-z'_i|^2+2C^2|y^-|^2.\ea\ee It
is easy to check that $\bar{\psi}^n(t,y,z)$ is also Lipschitz continuous in
$(y,z)$. From Theorem \ref{t4}, we know that RBSDE~(\ref{RBSDE9}) has a
unique adapted solution $(Y^n, Z^n, K^{-,n})$ with $$Y^n\triangleq
(Y^n_1,\cdots,Y^n_m)^T\in (\cS^2)^m,\q Z^n\triangleq
(Z^n_1,\cdots,Z^n_m)^T\in (\cM^2)^{m\times d},
$$
and
$$ K^{-,n}\triangleq (K^{-,n}_1,\cdots,K^{-,n}_m)^T\in
(\cN^{\,2})^m.$$ Moreover, we have from (\ref{89}) that for $(y,y')\in
(R^m)^2, (z,z')\in (R^{m\times d})^2,$ and $(i, l) \in (\L)^2,$
\be\ba{rcl}&&-4\langle
y^-,\bar{\psi}^{n+1}(t,y^++y',z)-\bar{\psi}^n(t,y',z')\rangle\\[3mm]&=&\displaystyle-4\langle
y^-,\bar{\psi}^n(t,y^++y',z)-\bar{\psi}^n(t,y',z')\rangle-4\sum_{i,\,l=1}^m\langle
y_i^-, (y_i^++y'_i-y_l^+-y'_l+k(i,l))^-\rangle\\[3mm]&\leq& \displaystyle-4\langle
y^-,\bar{\psi}^n(t,y^++y',z)-\bar{\psi}^n(t,y',z')\rangle\\&\leq&\displaystyle2\sum_{i=1}^m\chi_{\{y_i<0\}}|z_i-z'_i|^2+2C^2|y^-|^2.\ea\ee
From Theorem \ref{t7}, we know
 \be \label{monotone}Y^n(t)\leq Y^{n+1}(t),\q K^{-,n}(t)\leq K^{-,n+1}(t),\q 0\leq t\leq T.\ee For $t\in [0,T]$,
the sequence $\{(Y^n(t),K^{-,n}(t))\}_{n\geq1}$ admits a limit, which is
denoted by $(Y(t),K^-(t))$ below with $Y(t)\triangleq (Y_1(t), \cdots,
Y_m(t))^T$ and $K^-(t)\triangleq (K^-_1(t), \cdots, K^-_m(t))^T$.

\ms {\bf Step 2. A priori estimate.}

The following lemma is the key to our subsequent arguments.

 \bl\label{l7} There is a positive constant c which is
independent of n, such that \begin{eqnarray*}\ba{rcl}\displaystyle
E\sup_{0\leq t\leq T}|Y^n(t)|^2+E|K^{-,n}(T)|^2+E\int_0^T
|Z^n(s)|^2\, d s&\leq &
c,\\[3mm] \displaystyle
n^2 E\int_0^T|(Y_i^n(s)-Y_j^n(s)+k(i,j))^-|^2\, d s &\leq& c, \q
i\in\L.\ea\end{eqnarray*}\el

The proof follows. Applying It\^o-Meyer's  formula~\cite{Meyer} to compute
$|(Y_i^n(t)-Y_j^n(t)+k(i,j))^-|^2$,  we see that the sum  $
|(Y_i^n(t)-Y_j^n(t)+k(i,j))^-|^2+\int_t^T\chi_{L^-_{i,j,n}}(s)|Z_i^n(s)-Z_j^n(s)|^2\,
d s+2n\int_t^T|(Y_i^n(s)-Y_j^n(s)+k(i,j))^-|^2\, d s$ is equal to the
following \be\label{69}\ba{rcl}&&\displaystyle
2\int_t^T(Y_i^n(s)-Y_j^n(s)+k(i,j))^-[\psi(s,Y_j^n(s),Z_j^n(s),j)-\psi(s,Y_i^n(s),Z_i^n(s),i)]\,
d
s\\[3mm]&&\displaystyle+2\int_t^T(Y_i^n(s)-Y_j^n(s)+k(i,j))^-(Z_i^n(s)-Z_j^n(s))\, d
W(s)\\[3mm]&&\displaystyle+2\int_t^T(Y_i^n(s)-Y_j^n(s)+k(i,j))^-d(K_i^{-,n}(s)-K_j^{-,n}(s))\\&&\displaystyle+2n\int_t^T(Y_i^n(s)-Y_j^n(s)+k(i,j))^-(Y_j^n(s)-Y_i^n(s)+k(j,i))^-d
s\\[3mm]&&\displaystyle+2n\sum_{l\neq i,\, l\neq
j}\int_t^T(Y_i^n(s)-Y_j^n(s)+k(i,j))^-[(Y_j^n(s)-Y_l^n(s)+k(j,l))^-\\[3mm]&&\displaystyle-(Y_i^n(s)-Y_l^n(s)+k(i,l))^-]\, d
s\ea\ee where
$$L_{i,j,n}^-\triangleq \{(t,\o):
Y_i^n(t)-Y_j^n(t)+k(i,j)<0\}.$$

We claim that the last three terms of (\ref{69}) are all equal to or less
than $0$. In fact, due to (\ref{RBSDE9}), we have \be Y_j^n(s)-k(i,j)\le
S_j(s)-k(i,j)\le S_i(s) \q \hbox{\rm (noting $S(s)\in \tilde Q(s)$)}\ee and
\begin{eqnarray*}\ba{rcl}&&\displaystyle\int_t^T(Y_i^n(s)-Y_j^n(s)+k(i,j))^-d(K_i^{-,n}(s)-K_j^{-,n}(s))\\[3mm]&\leq&\displaystyle
\int_t^T(Y_i^n(s)-Y_j^n(s)+k(i,j))^-d
K_i^{-,n}(s)\\[3mm]&\leq&\displaystyle\int_t^T(Y_i^n(s)-S_i(s))^-d
K_i^{-,n}(s)\\[3mm]&=&\displaystyle\int_t^T(S_i(s)-Y_i^n(s))d K_i^{-,n}(s)=0.\ea\end{eqnarray*}  In view of Hypothesis 3'(i), we
have
$$\{(y_1,\cdots,y_m)^T\in R^m: y_i-y_j+k(i,j)<0,
y_j-y_i+k(j,i)<0\}=\emptyset,$$  which immediately gives  $$
(Y_i^n(t)-Y_j^n(t)+k(i,j))^-(Y_j^n(t)-Y_i^n(t)+k(j,i))^-=0. $$ From
Hypothesis 3'(ii), using the property that $x_1^--x_2^-\leq (x_1-x_2)^-$, we
have
\begin{eqnarray*}\ba{rcl}&&(Y_i^n(s)-Y_j^n(s)+k(i,j))^-[(Y_j^n(s)-Y_l^n(s)+k(j,l))^-\\[3mm]&&-(Y_i^n(s)-Y_l^n(s)+k(i,l))^-]\\[3mm]&\leq
&(Y_i^n(s)-Y_j^n(s)+k(i,j))^-(Y_j^n(s)-Y_i^n(s)+k(j,l)-k(i,l))^-\\[3mm]&\leq&(Y_i^n(s)-Y_j^n(s)+k(i,j))^-(Y_j^n(s)-Y_i^n(s)-k(i,j))^-\\[3mm]&=&(Y_i^n(s)-Y_j^n(s)+k(i,j))^-(Y_i^n(s)-Y_j^n(s)+k(i,j))^+=
0.\ea\end{eqnarray*}

Taking expectation on both sides of (\ref{69}), we have
\be\label{70}\ba{rcl}&&\displaystyle
E|(Y_i^n(t)-Y_j^n(t)+k(i,j))^-|^2+E\int_t^T\chi_{L^-_{i,j,n}}(s)|Z_i^n(s)-Z_j^n(s)|^2\,
d
s\\[3mm]&&\displaystyle+2n E\int_t^T|(Y_i^n(s)-Y_j^n(s)+k(i,j))^-|^2\, d
s\\[3mm]&\leq
&\displaystyle 2E\int_t^T
(Y_i^n(s)-Y_j^n(s)+k(i,j))^-|\psi(s,Y_i^n(s),Z_i^n(s),i)\\[3mm]&&-\psi(s,Y_j^n(s),Z_j^n(s),j)|\, d
s.\ea\ee Noting that
\begin{eqnarray*}\ba{rcl}&&|\psi(s,Y_i^n(s),Z_i^n(s),i)-\psi(s,Y_j^n(s),Z_j^n(s),j)|\\[3mm]&\leq&
|\psi(s,Y_i^n(s),Z_i^n(s),i)-\psi(s,Y_i^n(s),Z_i^n(s),j)|\\[3mm]&&+|\psi(s,Y_i^n(s),Z_i^n(s),j)-\psi(s,Y_j^n(s),Z_j^n(s),j)|\\[3mm]&\leq
&
c(|\psi(s,0,0)|+|Y_i^n(s)|+|Z_i^n(s)|+|Y_i^n(s)-Y_j^n(s)|+|Z_i^n(s)-Z_j^n(s)|)\\[3mm]&\leq&
c\(1+|\psi(s,0,0)|+|Y_i^n(s)|+|Z_i^n(s)|+|Y_i^n(s)-Y_j^n(s)+k(i,j)|\\&&+|Z_i^n(s)-Z_j^n(s)|\)\ea\end{eqnarray*}
for a positive constant $c$ (independent of $n$ and possibly varying
from line to line), in view of~(\ref{70}), we have
\begin{eqnarray*}\ba{rcl} &&\displaystyle
E|(Y_i^n(t)-Y_j^n(t)+k(i,j))^-|^2+E\int_t^T\chi_{L^-_{i,j,n}}(s)|Z_i^n(s)-Z_j^n(s)|^2d
s\\[3mm]&&\displaystyle+2n E\int_t^T|(Y_i^n(s)-Y_j^n(s)+k(i,j))^-|^2d
s\\[3mm]&\leq
&\displaystyle(\frac{n}{2}+2c)E\int_t^T|(Y_i^n(s)-Y_j^n(s)+k(i,j))^-|^2d
s+\frac{2c^2}{n}E\int_t^T\chi_{L^-_{i,j,n}}(s)\\[3mm]&&\displaystyle\(1+|\psi(s,0,0)|^2+|Y_i^n(s)|^2+|Z_i^n(s)|^2+|Z_i^n(s)-Z_j^n(s)|^2\)d
s.\ea\end{eqnarray*} So for sufficiently large $n$,
\be\label{73}\qq\q n^2 E\int_t^T|(Y_i^n(s)-Y_j^n(s)+k(i,j))^-|^2d s
\leq c\(1+E\int_t^T(|Y_i^n(s)|^2+|Z_i^n(s)|^2)d s\).\ee

Applying It\^{o}'s lemma to compute $|Y_i^n(t)|^2$, we have
\be\label{71}\ba{rcl}&&\displaystyle|Y_i^n(t)|^2+\int_t^T|Z_i^n(s)|^2d
s\\[3mm]&=&\displaystyle2\int_t^T
Y_i^n(s)\(\psi(s,Y_i^n(s),Z_i^n(s),i)+n\sum_{l=1}^m
(Y_i^n(s)-Y_l^n(s)+k(i,l))^-\)d s\\[3mm]&&\displaystyle+\xi_i^2-2\int_t^T Y_i^n(s)d
K_i^{-,n}(s)-2\int_t^T Y_i^n(s)Z_i^n(s)d W(s).\ea\ee Using the elementary
inequality:
$$\displaystyle2ab\leq
\frac{1}{\alpha}a^2+\alpha b^2 \q \hbox{ \rm for } a, b, \a>0,$$ we obtain
that \be\label{72}\ba{rcl}\displaystyle-2\int_t^T Y_i^n(s)\, d
K_i^{-,n}(s)&=&\displaystyle-2\int_t^T S_i(s)\, d K_i^{-,n}(s) \\[3mm]&\leq
&\displaystyle\frac{1}{\alpha}\sup_{0\leq t\leq T}|S_i^-(t)|^2+\alpha
|K_i^{-,n}(T)-K_i^{-,n}(t)|^2\ea\ee for  an arbitrary positive real number
$\alpha$. Then taking expectation on both sides of (\ref{71}), we have
\be\label{75}\ba{rcl} &&\displaystyle
E|Y_i^n(t)|^2+E\int_t^T|Z_i^n(s)|^2d s \\[3mm]&\leq &\displaystyle 2 E\int_t^T
|Y_i^n(s)|(|\psi(s,Y_i^n(s),Z_i^n(s),i)|+n\sum_{l=1}^m
(Y_i^n(s)-Y_j^n(s)+k(i,l))^-)d s\\&&\displaystyle+ E(\xi_i^2)-2 E\int_t^T
Y_i^n(s)d
K_i^{-,n}(s) \\[3mm]&\leq & \displaystyle c_\e
E\int_t^T|Y_i^n(s)|^2 d s+\e \sum_{l=1}^m n^2
E\int_t^T|(Y_i^n(s)-Y_l^n(s)+k(i,l))^-|^2 d
s+E(\xi_i)^2\\[3mm]&&\displaystyle+\e E\int_t^T|Z_i^n(s)|^2 d s+\frac{1}{\alpha}E\sup_{0\leq t\leq T}|S_i^-(t)|^2+\alpha
E|K_i^{-,n}(T)-K_i^{-,n}(t)|^2\ea\ee for arbitrary positive real numbers $\e,
\alpha,$ and  a constant $c_\e$ depending on $\e$.
 On the other hand, from
equation (\ref{RBSDE9}), we have \begin{eqnarray*}\ba{rcl}&&
K_i^{-,n}(T)-K_i^{-,n}(t)\\[3mm]&=&\displaystyle\xi_i-Y_i^n(t)+\int_t^T\psi(s,Y_i^n(s),Z_i^n(s),i)\, d
s\\[3mm]&&\displaystyle+n\sum_{l=1}^m \int_t^T(Y_i^n(s)-Y_l^n(s)+k(i,l))^-d s-\int_t^T
Z_i^n(s)\, d W(s).\ea\end{eqnarray*} In view of~(\ref{73}), we have
\be\label{74}\ba{rcl}&&E|K_i^{-,n}(T)-K_i^{-,n}(t)|^2\\[3mm]&\leq &\displaystyle
c\(E\xi_i^2+
E|Y_i^n(t)|^2+E\int_t^T(|Y_i^n(s)|^2+|Z_i^n(s)|^2)\,ds\)\\[3mm]&&\displaystyle+c\sum_{l=1}^m n^2E
\int_t^T|Y_i^n(s)-Y_l^n(s)+k(i,l))^-|^2ds \\[3mm]&\leq &\displaystyle
c \left[1+E|Y_i^n(t)|^2+E\int_t^T(|Y_i^n(s)|^2+|Z_i^n(s)|^2)\, d
s\right]. \ea\ee Further, in view of~(\ref{75}), we have
\begin{eqnarray*}\ba{rcl} &&\displaystyle
E|Y_i^n(t)|^2+E\int_t^T|Z_i^n(s)|^2d s \\[3mm]&\leq &\displaystyle c_{\e,\alpha}
(1+E\int_t^T|Y_i^n(s)|^2 d s)+\alpha c E|Y_i^n(t)|^2
+\frac{1}{\alpha}E\sup_{0\leq t\leq T}|S^-_i(t)|^2\\[3mm]&&\displaystyle+(\e+\alpha)c E\int_t^T |Z_i^n(s)|^2d s\ea\end{eqnarray*}
 for any positive real numbers $\e, \alpha$ and a constant $c_{\e,\alpha}$ depending on $\e, \alpha$.
 Setting
$\alpha=\e=\frac{1}{3c},$ we have from Gronwall's inequality that
$$E|Y_i^n(t)|^2+E\int_t^T|Z_i^n(s)|^2 d s\leq c.$$ From
(\ref{73}) and (\ref{74}), we have \be n^2
E\int_0^T|(Y_i^n(s)-Y_j^n(s)+k(i,j))^-|^2d s\leq c, \qq E|K_i^{+,n}(T)|^2\leq
c.\ee Moreover, in view of~(\ref{71}), applying the BDG inequality, we have
$$E\sup_{0\leq t\leq T}|Y_i^n(t)|^2\leq c.$$  The proof of Lemma~\ref{l7} is then
complete.

\ms {\bf Step 3. The convergence of penalized BSDEs.}

In view of the component-wisely monotone convergence of $\{(Y^n, K^{-,n}),
n=1,2,\ldots\}$ (see~(\ref{monotone})), we have
$$E[\sup_{0\leq t\leq T}(|Y(t)|^2+|K^-(t)|^2)]\le E[\mathop{\underline{\lim}}_{n\to \infty}\sup_{0\leq t\leq T}(|Y^n(t)|^2+|K^{-,n}(t)|^2)].$$
In view of Lemma~\ref{l7},  using Fatou's lemma, we have
$$E[\sup_{0\leq t\leq T}(|Y(t)|^2+|K^-(t)|^2)]\le \mathop{\underline{\lim}}_{n\to \infty} E[\sup_{0\leq t\leq T}(|Y^n(t)|^2+|K^{-,n}(t)|^2)]< \infty.$$ Then
applying Lebesgue's dominated convergence theorem, we have
\be\label{63}\qq E\int_0^T|Y^n(s)-Y(s)|^2\, ds + E\int_0^T
|K^{-,n}(s)- K^-(s)|^2\, d s\ra 0 \mbox{\q as \q } n\ra\infty.\ee
For positive integers $n_1$ and $n_2$, applying It\^{o}'s lemma to
$|Y_i^{n_1}(t)-Y_i^{n_2}(t)|^2$, we have
\be\label{61}\ba{rcl}&&\displaystyle|Y_i^{n_1}(t)-Y_i^{n_2}(t)|^2+\int_t^T|Z_i^{n_1}(s)-Z_i^{n_2}(s)|^2d
s\\[3mm]&=&\displaystyle2\int_t^T(\psi(s,Y_i^{n_1}(s),Z_i^{n_1}(s),i)-\psi(s,Y_i^{n_2}(s),Z_i^{n_2}(s),i))(Y_i^{n_1}(s)-Y_i^{n_2}(s))d
s\\[3mm]&&\displaystyle+2n_1\sum_{l=1}^m\int_t^T(Y_i^{n_1}(s)-Y_l^{n_1}(s)+k(i,l))^-(Y_i^{n_1}(s)-Y_i^{n_2}(s))ds\\[3mm]&&\displaystyle-2n_2\sum_{l=1}^m\int_t^T(Y_i^{n_2}(s)-Y_l^{n_2}(s)+k(i,l))^-(Y_i^{n_1}(s)-Y_i^{n_2}(s))ds\\[3mm]&&\displaystyle
-2\int_t^T(Y_i^{n_1}(s)-Y_i^{n_2}(s))d(K_i^{-,n_1}(s)-K_i^{-,n_2}(s))\\[3mm]&&\displaystyle-2\int_t^T(Y_i^{n_1}(s)-Y_i^{n_2}(s))(Z_i^{n_1}(s)-Z_i^{n_2}(s))d
W(s),\qq\forall \,i\in\L.\ea\ee Since \be\label{65}\ba{rcl} &&\displaystyle
\int_t^T(Y_i^{n_1}(s)-Y_i^{n_2}(s))d(K_i^{-,n_1}(s)-K_i^{-,n_2}(s))\\[3mm]&=&\displaystyle\int_t^T(S_i(s)-Y_i^{n_2}(s))d
K_i^{-,n_1}(s)+\int_t^T(S_i(s)-Y_i^{n_1}(s))dK_i^{-,n_2}(s)\geq 0,\q
i\in\L,\ea\ee in view of~(\ref{61}), we have for
$i\in\L$,\begin{eqnarray*}\ba{rcl}&&\displaystyle
E|Y_i^{n_1}(t)-Y_i^{n_2}(t)|^2+E\int_t^T|Z_i^{n_1}(s)-Z_i^{n_2}(s)|^2d
s\\[3mm]&=&\displaystyle2E\int_t^T(\psi(s,Y_i^{n_1}(s),Z_i^{n_1}(s),i)-\psi(s,Y_i^{n_2}(s),Z_i^{n_2}(s),i))(Y_i^{n_1}(s)-Y_i^{n_2}(s))d
s\\[3mm]&&\displaystyle+2n_1\sum_{l=1}^m E\int_t^T(Y_i^{n_1}(s)-Y_l^{n_1}(s)
+k(i,l))^-(Y_i^{n_1}(s)-Y_i^{n_2}(s))ds\\[3mm]&&\displaystyle-2n_2\sum_{l=1}^m E\int_t^T(Y_i^{n_2}(s)-Y_l^{n_2}(s)+k(i,l))^-(Y_i^{n_1}(s)-Y_i^{n_2}(s))ds
\\[3mm]\ea\end{eqnarray*}
which implies the following \be\label{62}\ba{rcl}&&\displaystyle
E|Y_i^{n_1}(t)-Y_i^{n_2}(t)|^2+E\int_t^T|Z_i^{n_1}(s)-Z_i^{n_2}(s)|^2d
s\\[3mm]
&\leq&\displaystyle 2C(C+1) E\int_t^T|Y_i^{n_1}(s)-Y_i^{n_2}(s)|^2d s
+\frac{1}{2}E\int_t^T|Z_i^{n_1}(s)-Z_i^{n_2}(s)|^2 d
s\\[3mm]

&&\displaystyle+2(E\int_0^T|Y_i^{n_1}(s)-Y_i^{n_2}(s)|^2d
s)^\frac{1}{2}\\[3mm]
&&\displaystyle\times  \left[\sum_{l=1}^m\sum_{h=1}^2n_h(E\int_0^T
|(Y_i^{n_h}(s)-Y_l^{n_h}(s)+k(i, l))^-|^2d s)^\frac{1}{2}\right].\ea\ee
Setting $t=0$,  in view of (\ref{63}) and Lemma \ref{l7}, we have
\be\label{68}E\int_0^T |Z_i^{n_1}(s)-Z_i^{n_2}(s)|^2d s\ra 0,\q \forall
\,i\in \L,\mbox{\q as \q } n_1,n_2\ra\infty.\ee So there exists $Z\triangleq
(Z_1,\cdots,Z_m)^T\in (\cM^2)^{m\times d}$ such that
$$ \lim_{n\ra\infty} E\int_0^T|Z_i^n(s)-Z_i(s)|^2d s=0,\q \forall \,i\in\L.$$
In view of (\ref{61}), applying the BDG inequality, we have
\be\label{67}E(\sup_{0\leq t\leq T}|Y_i^{n_1}(t)-Y_i^{n_2}(t)|^2)\ra
0,\q\forall \,i\in \L\mbox{\q as \q } n_1,n_2\ra\infty.\ee From
(\ref{RBSDE9}), we have \begin{eqnarray*}\ba{rcl}&&
d(K_i^{-,n_1}(t)-K_i^{-,n_2}(t))\\[3mm]&=&(\psi(t,Y_i^{n_1}(t),Z_i^{n_1}(t),i)-\psi(t,Y_i^{n_2}(t),Z_i^{n_2}(t),i))dt\\[3mm]
&&\displaystyle+d(Y_i^{n_1}(t)-Y_i^{n_2}(t))+n_1\sum_{l=1}^m(Y_i^{n_1}(t)-Y_l^{n_1}(t)+k(i,l))^-dt
\\[3mm]&&\displaystyle-n_2\sum_{l=1}^m(Y_i^{n_2}(t)-Y_l^{n_2}(t)+k(i,l))^-dt-(Z_i^{n_1}(t)-Z_i^{n_2}(t))d
W(t).\ea\end{eqnarray*} Since the process
$\{(K_i^{-,n_1}(t)-K_i^{-,n_2}(t)), t\in [0,T]\}$ is of finite
variation, its quadratic variation is 0. By It\^{o}'s lemma it
follows that
\begin{eqnarray*}\ba{rcl}&&(K_i^{-,n_1}(t)-K_i^{-,n_2}(t))^2\\[3mm]&=&\displaystyle 2\int_0^t(K_i^{-,n_1}(s)-K_i^{-,n_2}(s))(\psi(s,Y_i^{n_1}(s),Z_i^{n_1}(s),i)-\psi(s,Y_i^{n_2}(s),Z_i^{n_2}(s),i))d
s\\[3mm]&&\displaystyle+2\int_0^t(K_i^{-,n_1}(s)-K_i^{-,n_2}(s))(n_1\sum_{l=1}^m(Y_i^{n_1}(s)-Y_l^{n_1}(s)+k(i,l))^-)d
s\\[3mm]&&\displaystyle-2\int_0^t(K_i^{-,n_1}(s)-K_i^{-,n_2}(s))(n_2\sum_{l=1}^m(Y_i^{n_2}(s)-Y_l^{n_2}(s)+k(i,l))^-)d
s\\[3mm]&&\displaystyle+2\int_0^t(K_i^{-,n_1}(s)-K_i^{-,n_2}(s))d(Y_i^{n_1}(s)-Y_i^{n_2}(s))\\[3mm]&&\displaystyle-2\int_0^t(K_i^{-,n_1}(s)-K_i^{-,n_2}(s))(Z_i^{-,n_1}(s)-Z_i^{-,n_2}(s))d
W(s).\ea\end{eqnarray*}  Applying the BDG inequality, we have
\be\label{64}\ba{rcl}&&\displaystyle E[\sup_{0\leq t\leq
T}|K_i^{-, n_1}(t)-K_i^{-, n_2}(t)|^2]\\[3mm]&\leq&\displaystyle
c [E\int_0^T(|K_i^{-, n_1}(s)-K_i^{-,
n_2}(s)|^2+|Y_i^{n_1}(s)-Y_i^{n_2}(s)|^2+|Z_i^{n_1}(s)-Z_i^{n_2}(s)|^2)d
s]\\[3mm]
&&\displaystyle+2(E\int_0^T|K_i^{-, n_1}(s)-K_i^{-, n_2}(s)|^2d
s)^\frac{1}{2}\\[3mm]
&&\displaystyle\times  \left[\sum_{l=1}^m\sum_{h=1}^2n_h(E\int_0^T
|(Y_i^{n_h}(s)-Y_l^{n_h}(s)+k(i, l))^-|^2d
s)^\frac{1}{2}\right]\\[3mm]
&&\displaystyle+E\left[\sup_{0\leq t\leq
T}[2\int_0^t(K_i^{-, n_1}(s)-K_i^{-, n_2}(s))d(Y_i^{n_1}(s)-Y_i^{n_2}(s))]\right]\\[3mm]&&\displaystyle+\frac{1}{3}E[\sup_{0\leq
t\leq T}|K_i^{-, n_1}(t)-K_i^{-, n_2}(t)|^2],\ea\ee for a positive constant
$c$ independent of $n_1$ and $n_2$.   Identical to the proof of (\ref{65}),
we have
$$\int_0^t(Y_i^{n_1}(s)-Y_i^{n_2}(s))\, d(K_i^{-,n_1}(s)-K_i^{-,n_2}(s))\geq
0.$$ Hence by It\^{o}'s lemma,
\begin{eqnarray*}\ba{rcl}&&\displaystyle
2\int_0^t(K_i^{-,n_1}(s)-K_i^{-,n_2}(s))d(Y_i^{n_1}(s)-Y_i^{n_2}(s))\\[3mm]&=&
2
((K_i^{-,n_1}(s)-K_i^{-,n_2}(s))(Y_i^{n_1}(s)-Y_i^{n_2}(s)))\|_0^t\\[3mm]&&\displaystyle
-2\int_0^t(Y_i^{n_1}(s)-Y_i^{n_2}(s))d(K_i^{-,n_1}(s)-K_i^{-,n_2}(s))\\[3mm]&\leq&\displaystyle
2|K_i^{-,n_1}(t)-K_i^{-,n_2}(t)||Y_i^{n_1}(t)-Y_i^{n_2}(t)|.\ea\end{eqnarray*}
Then as a consequence, \be\label{66}\ba{rcl}&&\displaystyle
E\sup_{0\leq t\leq
T}\left\{2\int_0^t(K_i^{-,n_1}(s)-K_i^{-,n_2}(s))d(Y_i^{n_1}(s)-Y_i^{n_2}(s))\right\}\\[3mm]&\leq
& \displaystyle 2 E\(\sup_{0\leq t\leq
T}|K_i^{-,n_1}(t)-K_i^{-,n_2}(t)||Y_i^{n_1}(t)-Y_i^{n_2}(t)|\)\\[3mm]&\leq&\displaystyle
\frac{1}{3}E\sup_{0\leq t\leq T}|K_i^{-,n_1}(t)-K_i^{-,n_2}(t)|^2+3E
\sup_{0\leq t\leq T}|Y_i^{n_1}(t)-Y_i^{n_2}(t)|^2.\ea\ee Together
with (\ref{63}), (\ref{64}) and Lemma~\ref{l7}, we have
\begin{eqnarray*}\ba{rcl}&&\displaystyle E[\sup_{0\leq
t\leq T}|K_i^{-, n_1}(t)-K_i^{-, n_2}(t)|^2]\\[3mm]&\leq & \displaystyle
c E[\sup_{0\leq t\leq
T}|Y_i^{n_1}(t)-Y_i^{n_2}(t)|^2+\int_0^T(|Y_i^{n_1}(s)-Y_i^{n_2}(s)|^2+|Z_i^{n_1}(s)-Z_i^{n_2}(s)|^2)d
s]. \ea\end{eqnarray*} From (\ref{68}) and (\ref{67}), we have
\begin{eqnarray*} E\sup_{0\leq t\leq
T}|K_i^{-,n_1}(t)-K_i^{-,n_2}(t)|^2\ra 0, \q \forall i\in\L,\mbox{\q
as\q} n_1,n_2\ra\infty.\end{eqnarray*}
 Set
$$\displaystyle K_i^{+,n}(t)\triangleq
n\sum_{l=1}^m\int_0^t(Y_i^n(s)-Y_l^n(s)+k(i,l))^- d s$$ and
\begin{eqnarray*} K_i^+(t)\triangleq
Y_i(0)-Y_i(t)-\int_0^t\psi(s,Y_i(s),Z_i(s),i)d s+\int_0^t d
K_i^-(s)+\int_0^t Z_i(s)d W(s).\end{eqnarray*} We have
\begin{eqnarray*}
K_i^{+,n}(t)=Y_i^n(0)-Y_i^n(t)-\int_0^t\psi(s,Y_i^n(s),Z_i^n(s),i)d
s+\int_0^t d K_i^{-,n}(s)+\int_0^tZ_i^n(s)d W(s)\end{eqnarray*} and
\begin{eqnarray*}\ba{rcl}&&\displaystyle \lim_{n\ra\i}E(\sup_{0\leq t\leq
T}|K_i^{+,n}(t)-K_i^+(t)|^2)\\[3mm]&\leq&\displaystyle \lim_{n\ra\i}c E\(\sup_{0\leq t\leq
T}\{|Y_i^n(t)-Y_i(t)|^2+|K_i^{-,n}(t)-K_i^-(t)|^2\}\\[3mm]&&\displaystyle+\int_0^T|Z_i^n(s)-Z_i(s)|^2
d s
\)\\[3mm]&=&0.\ea\end{eqnarray*}
Hence, $K^+\triangleq (K_1^+,\cdots,K_m^+)^T\in(\cN^{\,2})^m$, and
$(Y,Z,K^+,K^-)$ satisfies the first equation of (\ref{DRBSDE8}).
From Lemma \ref{l7}, using Fatou's lemma, we obtain that
\begin{eqnarray*}\ba{rcl} &&\displaystyle
E\int_0^T|(Y_i(s)-Y_l(s)+k(i,l))^-|^2d
s\\[3mm]&\leq& \displaystyle
\mathop{\underline{\lim}}_{n\ra\infty}E\int_0^T|(Y_i^n(s)-Y_l^n(s)+k(i,l))^-|^2d s\\[3mm]&\leq
&\displaystyle\lim_{n\ra\infty}\frac{c}{n^2}=0,\ea\end{eqnarray*}
which implies immediately that \be\label{99}Y_i(t)-Y_l(t)+k(i,l)\geq
0, \q 0\leq t\leq T.\ee Since $Y^n(t)\leq S(t)$ for any $n$ and
$t\in [0,T]$,  we have
$$Y(t)\leq S(t), \q 0\leq t\leq T.$$ Hence,  $$Y(t)\in \tilde Q(t),  \q
0\leq t\leq T.$$

\ms {\bf Step 4. The lower minimal boundary condition.}

Consider $(i, j, l)\in\L\times\L\times\L$ such that $i\neq j$. When
$l=i$, it is obvious that \begin{eqnarray*} (Y_i^n(s)-Y_j^n(s)+k(i,
j))^+(Y_i^n(s)-Y_l^n(s)+k(i, l))^-=0. \end{eqnarray*}  If $l\neq i$,
we have \begin{eqnarray*}\ba{rcl}&&\min\limits_{j\neq
i}\{(Y_i^n(s)-Y_j^n(s)+k(i, j))^+(Y_i^n(s)-Y_l^n(s)+k(i, l))^-\}\\[3mm]&\leq &
(Y_i^n(s)-Y_l^n(s)+k(i, l))^+(Y_i^n(s)-Y_l^n(s)+k(i, l))^-=0.
\ea\end{eqnarray*} Hence,
\begin{eqnarray*}\ba{rcl}&&\displaystyle\int_0^T
(Y_i^n(s)-\max_{j\neq i}\{Y_j^n(s)-k(i,j)\})^+d
K_i^{+,n}(s)\\[3mm]&=&\displaystyle
n\sum_{l=1}^m \int_0^T\min_{j\neq
i}\{(Y_i^n(s)-Y_j^n(s)+k(i,j))^+(Y_i^n(s)-Y_l^n(s)+k(i,j))^-\} d
s\\[3mm]
&\leq& 0.\ea\end{eqnarray*} On the other hand, since
$K_i^{+,n}(\cdot)$ is increasing,  we have
$$\int_0^T (Y_i^n(s)-\max_{j\neq i}\{Y_j^n(s)-k(i,j)\})^+d
K_i^{+,n}(s)\geq 0.$$ Therefore,
$$\int_0^T (Y_i^n(s)-\max_{j\neq i}\{Y_j^n(s)-k(i,j)\})^+d
K_i^{+,n}(s)=0.$$ Following the same arguments to the proof of (\ref{90}) and
applying \cite[Lemma 5.8]{Gegout-Petit}, we have
$$\int_0^T (Y_i(s)-\max_{j\neq i}\{Y_j(s)-k(i,j)\})^+d K_i^+(s)=0.$$
In view of  (\ref{99}), we have $$\int_0^T (Y_i(s)-\max_{j\neq
i}\{Y_j(s)-k(i,j)\})d K_i^+(s)=0.$$
 \qquad\end{proof}

\section{Uniqueness}
Uniqueness of the solution of RBSDE~(\ref{DRBSDE}) is defined in the
following sense:

If $(Y',Z',K'^{+},K'^{-})$ is another solution, then $$Y'(t)\equiv Y(t),\,
Z'(t)\equiv Z(t),\, K'^{+}(t)-K'^{-}(t)\equiv K^{+}(t)-K^{-}(t),\, \forall
\,0\leq t\leq T, a.s.$$

The following stronger assumption on $h_{i,j}$ is needed in our proof of the
uniqueness result.

 \textbf{Hypothesis 5.} For $(i,j,l)\in\L\times\L\times\L$ such that $i\neq j,\, j\neq l$,
$$h_{i,j}(t,h_{j,l}(t,y))< h_{i,l}(t,y),\q \forall \,y\in
R.$$

\br It is easy to check that Hypothesis 5 implies Hypothesis 3. If
$h_{i,j}(t,y)\triangleq y-k(i,j)$, then Hypothesis 5 reduces to the
inequality: $k(i,j)+k(j,l)>k(i,l)$ for  $i\neq j,\,j\neq l$.\er

 Let $\{\theta_j\}_{j=0}^{\i}$ be an increasing sequence of
stopping times with values in $[0,T]$. $\forall \,j,\, \alpha_j$ is an
$\cF_{\theta_j}$-measurable random variable with values in $\L$. Assume that
$a.s. \,\o$, there exists an integer ${\wt N}(\o)<\i$ such that
$\theta_{\tilde N}=T$. Define $N(\o)$ to be the smallest integer ${\wt
N}(\o)<\i$ such that $\theta_{\tilde N}=T$. Then the set $\{N=j\}$ is
$\cF_{\theta_j}$-measurable for any integer $j$.  Then we define an
admissible switching control as:
$$a(s)=\alpha_0\chi_{[\theta_0,\,\theta_1]}(s)+\sum_{j=1}^{N-1}\alpha_j\chi_{(\theta_j,\, \theta_{j+1}]}(s),
\q s\in[\theta_0,\,T].$$ Denote  by $\cA_t^i$ all the admissible switching
controls with initial data $(\alpha_0, \theta_0)=(i,t)\in\L \times [0,T]$,
and by $\cT_t$ the totality of stopping times which take values in $[t,T]$.
For given $a\in\cA_t^i$, consider RBSDE~(\ref{BSDE1}).
 The generator $\psi$ of RBSDE~(\ref{BSDE1}) depends
on the control $a$, and at each switching time $\theta_j$ before termination,
the value of $U^a$ will jump by an amount of
$U^a(\theta_j)-h_{\alpha_{j-1},\,\alpha_j}(\theta_j,U^a(\theta_j)),$ which
can be regarded as a penalty or cost for the switching. In each subinterval
divided by the switching times, RBSDE (\ref{BSDE1}) evolves as a standard
RBSDE with single barrier, which can be solved in a backwardly inductive way.

  The optimal control
problem for RBSDE~(\ref{BSDE1}) is  to maximize $U^a(t)$ over $a\in \cA_t^i$.
The solution of RBSDE~(\ref{DRBSDE}) is closely connected with this control
problem. Besides, it is also connected to the stochastic game constructed
below.

 For given $a\in\cA_t^i$ and stopping time $\tau\in \cT_t$, consider the following BSDE:
\be\label{BSDE5.2}\ba{rcl} U^{a,\tau}(s)&=&\displaystyle
S_{a(\tau)}(\tau)\chi_{\{\tau<T\}}+\xi_{a(\tau)}\chi_{\{\tau=T\}}+\int_s^{\tau}\psi(r,U^{a,\tau}(r),V^{a,\tau}(r),a(r))d
r\\[3mm]&&\displaystyle-\sum\limits_{j=1}^{N-1}
 [U^{a,\tau}(\theta_j)-h_{\alpha_{j-1},\,\alpha_j}(\theta_j,U^{a,\tau}(\theta_j))
]\chi_{(s,\tau]}(\theta_j)\\[3mm]&&\displaystyle-\int_s^{\tau}V^{a,\tau}(r)d
W(r),\qq s\in [t,\tau].\ea\ee Since the terminal value
$S_{a(\tau)}(\tau)\chi_{\{\tau<T\}}+\xi_{a(\tau)}\chi_{\{\tau=T\}}$  depends
on the terminal time $\tau$, the first solution component $U^{a,\tau}(\cdot)$
depends on both the switching strategy $a$ and terminal time $\tau$. Like
RBSDE~(\ref{BSDE1}),  BSDE~(\ref{BSDE5.2}) can be shown in a backwardly
recursive way to  have a unique solution. We construct a zero-sum stochastic
game as follows. Suppose that there are two players A and B whose benefits
are antagonistic. The payoff $U^{a,\tau}(t)$ which is defined by
BSDE~(\ref{BSDE5.2}) is a reward for player A and a cost for player B. Player
A has the right to choose a switching strategy $a\in \cA_t^i$ so as to
maximize the reward $U^{a,\tau}(t)$. Player B has the right to choose the
time $\tau$ to terminate the game and tries to minimize the cost
$U^{a,\tau}(t)$.

The following theorem reveals the connections among the solution of the
system (\ref{DRBSDE}), the above control problem and stochastic game.

\bt\label{t3} Let Hypotheses 1 and 5 be satisfied. Assume that $(U^a, V^a,
L^a)$ is the unique adapted solution of RBSDE~(\ref{BSDE1}) for $a\in
\cA_t^i$ and $(U^{a,\tau}, V^{a,\tau})$ is the unique solution of
BSDE~(\ref{BSDE5.2}) for $a\in \cA_t^i$ and stopping time $\tau$. Then if
$(Y,Z,K^+,K^-)\in (\cS^2)^m\times (\cM^2)^{m\times d}\times (\cN^{\,
2})^{2m}$ is an adapted  solution of RBSDE~(\ref{DRBSDE}), we have
$$Y_i(t)=\esssup_{a\in\cA_t^i}U^a(t)=\esssup_{a\in\cA_t^i}\essinf_{\tau\in \cT_t} U^{a,\tau}(t)
=\essinf_{\tau\in \cT_t}\esssup_{a\in\cA_t^i}U^{a,\tau}(t),\q t\in[0,T].$$

Let $\hat{\theta}_0=t$ and $\hat{\alpha}_0=i$.  Define the sequence
$\{\hat{\theta}_j,\hat{\alpha}_j\}_{j=1}^{\i}$ in an inductive way as
follows: $$\hat{\theta}_j=\inf\{s\geq \hat{\theta}_{j-1}:
Y_{\hat{\alpha}_{j-1}}(s)=\max_{k\neq \hat{\alpha}_{j-1},k\in\L}
h_{\hat{\alpha}_{j-1},k}(s,Y_k(s))\}\wedge T.$$ And if $\hat{\theta}_j<T$,
set $\hat{\alpha}_j$  be the smallest index in $\L$ such that
\be\label{21}Y_{\hat{\alpha}_{j-1}}(\hat{\theta}_j)=
h_{\hat{\alpha}_{j-1},\hat{\alpha}_j}(\hat{\theta}_j,Y_{\hat{\alpha}_j}(\hat{\theta}_j)).\ee
Otherwise, set $\hat{\alpha}_j$ be an arbitrary index. Define
$$\hat{a}(s)\triangleq\hat{\alpha}_0\chi_{[\hat{\theta}_0,\,\hat{\theta}_1]}(s)+\sum_{j=1}^{\hat{N}-1}\hat{\alpha}_j\chi_{(\hat{\theta}_j,\,
\hat{\theta}_{j+1}]}(s)$$ and
$$\tau^*\triangleq \inf\{s\in[t,T): U^{\hat a}(s)=S_{{\hat a}(s)}(s)\}\wedge T,$$ with the convention that $\inf\emptyset\triangleq
+\i.$ Then, we have $\hat a\in \cA_t^i$ and \be Y_i(t)=U^{\hat a}(t)\q \hbox{
\rm and } \q Y_i(t)=U^{\hat a, \tau^*}(t), \q t\in [0,T]. \ee
 \et

\begin{proof} For $t\in[0,T]$ and $a\in \cA_t^i$, define \be\label{5.1}\ba{rcl}
Y^a(s)&\triangleq&\displaystyle\sum_{j=1}^N
Y_{\alpha_{j-1}}(s)\chi_{[\theta_{j-1},\,
\theta_j)}(s)+\xi_{\alpha_{N-1}}\chi_{\{s=T\}},\q s\in[t,T]; \\[3mm]
Z^a(s)&\triangleq&\displaystyle\sum_{j=1}^N
Z_{\alpha_{j-1}}(s)\chi_{[\theta_{j-1},\, \theta_j)}(s), \q s\in[t,T];\\[3mm]
K^{+,a}(s)&\triangleq&\displaystyle\sum_{j=1}^N
(K_{\alpha_{j-1}}^+(\theta_j\wedge s)-K_{\alpha_{j-1}}^+(\theta_{j-1}\wedge
s)), \q s\in[t,T];\\[3mm]K^{-,a}(s)&\triangleq&\displaystyle\sum_{j=1}^N (K_{\alpha_{j-1}}^-(\theta_j\wedge
s)-K_{\alpha_{j-1}}^-(\theta_{j-1}\wedge s)), \q s\in[t,T].\ea\ee In view of
the jump
$Y_{\alpha_j}(\theta_j)-Y_{\alpha_{j-1}}(\theta_j)=Y^a(\theta_j)-Y^a(\theta_j-)$
at each stopping time $\theta_j, \, j=1,\cdots,N-1$, we know that $(Y^a, Z^a,
K^{+,a}, K^{-,a})$ satisfies the following RBSDE:
\be\label{BSDE4}\left\{\ba{rcl}
Y^a(s)&=&\displaystyle \xi_{a(T)}+\int_s^T\psi(r,Y^a(r),Z^a(r),a(r))d r\\[3mm]&&\displaystyle-\sum\limits_{j=1}^{N-1}(Y^a(\theta_j)-Y^a(\theta_j-))\chi_{(s,T]}(\theta_j)+ K^{+,a}(T)-K^{+,a}(s)\\[3mm]&&\displaystyle
-(K^{-,a}(T)-K^{-,a}(s))-\int_s^T Z^a(r)d
W(r),\qq s\in [0,T]; \\[3mm]
Y^a(s)&\leq & S_{a(s)}(s), \qq s\in [0,T];\\[3mm]&&\displaystyle\int_t^T (Y^a(s)-S_{a(s)}(s))d K^{-,a}(s)=0,\qq s\in[t,T].\ea\right.\ee
For any stopping time $\tau\in \cT_t$, $(Y^a, Z^a, K^{+,a}, K^{-,a})$ also
satisfies the following BSDE: \be\label{BSDE4.1}\ba{rcl}
Y^a(s)&=&\displaystyle Y^a(\tau)+\int_s^{\tau}\psi(r,Y^a(r),Z^a(r),a(r))d r\\[3mm]&&\displaystyle
-\sum\limits_{j=1}^{N-1}(Y^a(\theta_j)-Y^a(\theta_j-))\chi_{(s,\tau]}(\theta_j)+K^{+,a}(\tau)-K^{+,a}(s)\\[3mm]&&\displaystyle-(K^{-,a}(\tau)-K^{-,a}(s))-\int_s^{\tau}Z^a(r)d
W(r),\qq s\in[t,\tau].\ea\ee

Comparing RBSDE~(\ref{BSDE1}) with (\ref{BSDE4}), in view of the facts that
$$Y^a(\theta_j-)=Y_{\alpha_{j-1}}(\theta_j)\geq
h_{\alpha_{j-1},\alpha_j}(\theta_j,Y_{\alpha_j}(\theta_j))=h_{\alpha_{j-1},\alpha_j}(\theta_j,Y^a(\theta_j))$$
and $$K^{+,a}(T)-K^{+,a}(t)\geq 0,$$ we deduce from the comparison theorem
\cite[Theorem 4.1]{El Karoui} that
$$Y^a(s)\geq U^a(s),\q  t\leq s\leq T.$$ From the definition in (\ref{5.1}),
$$ Y^a(t)=Y_i(t).$$ Hence,
\be\label{5.21}Y_i(t)\geq U^a(t), \q \forall a\in\cA_t^i.\ee

For the sequence $\{\hat{\theta}_j\}_{j= 1}^{\i}$, we claim that for $a.s.\,
\o$, there exists an integer $\hat{N}(\o)<\i$ such that
$\hat{\theta}_{\hat{N}}=T$. Otherwise, define $B\triangleq
\bigcap_{j=1}^{\i}\{\o:\hat{\theta}_j(\o)<T\}\in \cF_T$, then $P(B)>0$. For
$j=1,2,\cdots$, we have \be\label{20}\ba{rcl}
Y_{\hat{\alpha}_{j-1}}(\hat{\theta}_j)&=&
h_{\hat{\alpha}_{j-1},\hat{\alpha}_j}(\hat{\theta}_j,Y_{\hat{\alpha}_j}(\hat{\theta}_j)),\\[3mm]Y_{\hat{\alpha}_j}(\hat{\theta}_{j+1})&=&
h_{\hat{\alpha}_j,\hat{\alpha}_{j+1}}(\hat{\theta}_{j+1},Y_{\hat{\alpha}_{j+1}}(\hat{\theta}_{j+1})),\qq
\mbox{\rm on }\, B.\ea\ee Since the sequence $\{(\hat{\alpha}_{j-1},
\hat{\alpha}_j, \hat{\alpha}_{j+1})\}_{j=1}^{\i}$  takes values in $\L^3$,
which is a finite set, there are a triplet $(i_1, i_2, i_3)$ and a
subsequence $j_k$ such that for $ k=1,2\cdots$,
$$ (\hat{\alpha}_{j_k-1},\hat{\alpha}_{j_k},
\hat{\alpha}_{j_k+1})=(i_1,i_2,i_3).$$
 Since the sequence
$\{\hat{\theta}_j\}_{j=1}^{\i}$ is increasing and bounded by T, there is a
limit $\hat{\theta}_{\i}$. Passing to the limit in~(\ref{20}) for the
subsequence $\{j_k\}$, we have \be\ba{rcl}Y_{i_1}(\hat{\theta}_{\i})&=&
h_{i_1,i_2}(\hat{\theta}_{\i},Y_{i_2}(\hat{\theta}_{\i})),\\[3mm]Y_{i_2}(\hat{\theta}_{\i})&=&
h_{i_2,i_3}(\hat{\theta}_{\i},Y_{i_3}(\hat{\theta}_{\i})),\qq\qq\mbox{\rm
on }\,B.\ea\ee From Hypothesis 5, we have  \begin{eqnarray*}\ba{rcl}
Y_{i_1}(\hat{\theta}_{\i})&=&
h_{i_1,i_2}(\hat{\theta}_{\i},h_{i_2,i_3}(\hat{\theta}_{\i},Y_{i_3}(\hat{\theta}_{\i})))\\[3mm]&<
&h_{i_1,i_3}(\hat{\theta}_{\i},Y_{i_3}(\hat{\theta}_{\i}))\qq\qq
\mbox{\rm on }\,B.\ea\end{eqnarray*} This contradicts to the fact
that \be Y_{i_1}(\hat \theta_\infty)\geq\max\limits_{j\neq i_1, j\in
\L}h_{i_1,j}(\hat\theta_\infty,Y_j(\hat\theta_\infty)).\ee This
shows $\hat a \in \cA_t^i.$

It is easy to see that \be\label{5.22}Y^{\hat{a}}(\hat{\theta}_j-)=
h_{\hat{\alpha}_{j-1},\hat{\alpha}_j}(\hat{\theta}_j,Y^{\hat{a}}(\hat{\theta}_j)),\q
 K^{+,\hat{a}}(T)-K^{+,\hat{a}}(t)=0.\ee  Then $(Y^{\hat{a}}, Z^{\hat{a}},
K^{-,\hat{a}})$ satisfies RBSDE~(\ref{BSDE1}). By the uniqueness of the
solution of RBSDE~(\ref{BSDE1}), we have \be\label{100}
Y^{\hat{a}}(s)=U^{\hat{a}}(s),\q t\leq s\leq T.\ee Hence,
$$Y_i(t)=U^{\hat{a}}(t).$$  Noting (\ref{5.21}), we know
\be\label{5.24}Y_i(t)=\esssup_{a\in\cA_t^i}U^a(t),\qq \forall t\in[0,T].\ee

Define
$$\hat{\tau}\triangleq \inf\{s\in[t,T): Y^{\hat a}(s)=S_{{\hat a}(s)}(s)\}\wedge T,$$ with the convention that $\inf\emptyset\triangleq
+\i.$ From the upper minimal boundary condition in~(\ref{DRBSDE}), we have
\be\label{19}\ba{l}K^{-,{\hat a}}(\hat{\tau})-K^{-,{\hat a}}(s)=0,\\[3mm]
Y^{\hat a}(\hat{\tau})=S_{{\hat a}(\hat
\tau)}(\hat{\tau})\chi_{\{\hat{\tau}<T\}}+\xi_{{\hat
a}(\hat{\tau})}\chi_{\{\hat{\tau}=T\}}.\ea\ee
 In view
of~(\ref{BSDE4.1}) and~(\ref{5.22}), we have \be\label{BSDE7}\ba{rcl}
Y^{\hat{a}}(s)&=&\displaystyle S_{\hat
a(\hat\tau)}(\hat{\tau})\chi_{\{\hat{\tau}<T\}}+\xi_{\hat{a}(\hat{\tau})}\chi_{\{\hat{\tau}=T\}}+\int_s^{\hat{\tau}}\psi(r,Y^{\hat{a}}(r),Z^{\hat{a}}(r),{\hat{a}}(r))d
r\\[3mm]&&\displaystyle
-\sum\limits_{j=1}^{\hat{N}}(Y^{\hat{a}}(\hat{\theta}_j)-h_{\hat{\alpha}_{j-1},\hat{\alpha}_j}(\hat{\theta}_j,Y^{\hat{a}}(\hat{\theta}_j)))\chi_{(s,\hat{\tau}]}(\hat{\theta}_j)\\
[3mm]&&\displaystyle -\int_s^{\hat\tau}Z^{\hat{a}}(r)\, d W(r), \q
s\in [t, \hat \tau].\ea\ee Comparing the last equation with
BSDE~(\ref{BSDE5.2}), in view of the facts that \be
Y^{\hat{a}}(\tau)\leq S_{\hat
a(\tau)}(\tau)\chi_{\{\tau<T\}}+\xi_{\hat{a}(\tau)}\chi_{\{\tau=T\}}\ee
and $K^{+,a}$ and that $K^{-,\hat{a}}$ are increasing processes, we
have \begin{eqnarray*} Y^{\hat{a}}(t)\leq U^{\hat{a},\tau}(t),\q
Y^{\hat{a}}(t)= U^{\hat{a},\hat{\tau}}(t),\q
t\in[0,T].\end{eqnarray*} Noting by definition that
$$Y^{\hat{a}}(t)=Y_i(t),\q  t\in[0,T],$$ we have \be\label{25}
Y_i(t)=U^{\hat{a},\hat{\tau}}(t),\q  U^{\hat{a},\hat{\tau}}(t)\leq
U^{\hat{a},\tau}(t),\q  t\in [0,T].\ee Similar to the proof
of~(\ref{5.24}), we have \be Y_i(t)= \esssup_{a\in\cA_t^i}U^{a, \hat
\tau}(t),\qq \forall t\in[0,\hat \tau]\ee which implies \be
U^{a,\hat{\tau}}(t)\leq U^{\hat{a},\hat{\tau}}(t), \qq a\in
 \cA_t^i.\ee In view of~(\ref{25}), we conclude that
$(\hat{a},\hat{\tau})$ is a saddle point for the functional
$U^{a,\tau}(t), \\  (a,\tau)\in \cA_t^i\times \cT_t$:
\be\label{5.25} \qq
Y_i(t)=U^{\hat{a},\hat{\tau}}(t)=\esssup_{a\in\cA_t^i}\essinf_{\tau\in
\cT_t} U^{a,\tau}(t)=\essinf_{\tau\in
\cT_t}\esssup_{a\in\cA_t^i}U^{a,\tau}(t),\q \forall t\in[0,T].\ee
\qquad\end{proof}

Theorem \ref{t3} gives the uniqueness of the first solution component $Y$.
Uniqueness of other components $(Z,K^+,K^-)$ of the solution $(Y,Z,K^+,K^-)$
is a consequence of Doob-Meyer decomposition of $Y$. We conclude  the
following results.

\bt\label{t8} Let Hypotheses 1, 2 and 5 be satisfied. Assume that the upper
barrier $S$ is super-regular with $S(t)\in Q(t)$ for $t\in [0,T]$, and that
the terminal value $\xi\in L^2(\O,\cF_T,P;R^m)$ takes values in $Q(T)$. Then
RBSDE~(\ref{DRBSDE}) has a unique adapted solution $(Y,Z,K^+,K^-)\in
(\cS^2)^m\times (\cM^2)^{m\times d}\times (\cN^{\, 2})^{2m}$. \et

 \bt Let
Hypotheses 1 and 3'(i) be satisfied. Assume  that the upper barrier $S\in
(\cS^2)^m$  with $S(t)\in \tilde Q(t)$ for $t\in [0,T]$, that the terminal
value $\xi\in L^2(\O,\cF_T,P;R^m)$ takes values in $\tilde Q(T)$, and that
$k(i,j)+k(j,l)>k(i,l)$ for $i\neq j,\,j\neq l, \,i,j,l\in\L$. Then
RBSDE~(\ref{DRBSDE8}) has a unique adapted solution in $(\cS^2)^m\times
(\cM^2)^{m\times d}\times (\cN^{\, 2})^{2m}$. \et

\section{Economic interpretation as real options}

Switched RBSDEs can be interpreted as real options. Real options are real,
fully or partially  irreversible, investment opportunities which incur sunk
costs. In this sense they are similar to sequences of American-style
financial options with the sunk costs being regarded as strike prices, and
investors not only have flexibility to choose the time of investment, but
also have the flexibility to choose the scale or mode of investment. Brennan
and Schwartz \cite{Brennan} and MacDonald and Siegel \cite{McDonald} made
early contributions to real options analysis. Brennan and Schwartz applied
the idea of real options to valuation and exploration of natural resources,
and MacDonald and Siegel studied the value of waiting to invest. Dixit and
Pindyck \cite{Dixit} gives an excellent exposition on the topic.

The interpretation of switched RBSDEs as real options are as follows. An
investment opportunity is given to an economic agent, with the value being
equal to $U^a(S)$ which depends on the agent's operation process $a$. The set
$\Lambda$ is the set of all the  modes of operations for the investment
opportunity and $a(s)$, with values in $\Lambda$, is the mode of operation by
the agent at time $s$. The investment can be viewed as a complex option. The
agent has the right to choose both the time to change the mode of operation
and a new mode for the next operation.   Switching from $\alpha_{j-1}$ to
$\alpha_j$ can be viewed as exercising an American-type financial option with
the strike price  being the sunk cost
$U^a(\theta_j)-h_{\alpha_{j-1},\,\alpha_j}(\theta_j,U^a(\theta_j))$.
 In traditional
real option models the number of modes of operation is either $2$, invest or
do not invest, or $3$, do not enter, enter, and exit. The models in this
paper and Hu and Tang~\cite{Hu2}  are thus far-reaching generalizations of
traditional real option models to the case where there are $m$-different
modes of operations. An interesting example is the case of a staged
development of enterprise. Numbers in $\Lambda$ correspond to stages of
development in the venture enterprise and $U^a(s)$ is the value of the
enterprise in this example. When the enterprise is maturing at stage $j-1$
(e.g., product development), it is optimal for the enterprise to move into
the next stage, stage $j$, (e.g., marketing) by paying a sunk cost. The
models are flexible enough to allow to skip an intermediate stages if the
enterprise is successful enough or to go back to previous stages if the
economic environment deteriorates.

One notable feature of the model in this paper which distinguishes itself
from that in~Hu and Tang~\cite{Hu2} is the existence of the upper barrier,
$S_{a(s)}(s)$, for the value process, $U^a(s)$. An economic interpretation is
as follows. We may consider the model as a game between two players. One is
the economic agent who controls the mode of operations, and the other has the
right to terminate the project by paying an amount money equal to
$S(\tau)\chi_{\tau<T}+\xi\chi_{\tau=T}$. In the previous example of a staged
development of a venture enterprise, the second player may be a venture
capitalist who is willing to take over the company by paying the designated
amount at an optimal time to minimize the value of the opponent in the real
option. The second player's problem is mathematically an optimal stopping
problem to minimize his opponent's project value, and the associated value
process is known to be characterized by an RBSDE (see El Karoui et.
al.~\cite{El Karoui} for more details).

There is an extensive literature on games with real option features.
Interested readers are referred to Dixit and Pindyck~\cite[Chapter 9]{Dixit},
Bensoussan, Diltz, and Hoe~\cite{Bensoussan(a)},   Grenadier
\cite{Grenadier(a),Grenadier(b),Grenadier(c),Grenadier(d)}, and Smit and
Trigeorgis~\cite{Smit}.

\section*{Acknowledgment}
The authors thank the three referees for their helpful comments and suggestions.

\bibliographystyle{siam}

\end{document}